# CONSTRAINED LOCAL EXTREMA WITHOUT LAGRANGE MULTIPLIERS AND THE HIGHER DERIVATIVE TEST

SALVADOR GIGENA

*Dedicated to the future generations of Mathematicians. Mainly, and with very special affection, to Camila, Cristian, Pamela, Santino, Lautaro, Martina (Gabrielita) and Brenda.*

**Abstract.** – *Two are the main objectives of this article: first, we introduce a method for determining and analyzing constrained local extrema that provides a different alternative to all previous works on the topic, by eliminating Lagrange multipliers and reducing constrained problems to unconstrained ones. It has two added advantages: is very straightforward, effective and computationally faster than all of the previous ones; is also very simple in its conception and in the requirements of theoretical background material needed for its efficient application. In fact, to be performed successfully by undergraduate students, they only need a basic knowledge of three very central results from classical mathematical analysis: implicit function theorem, chain rule and Taylor's formula, together with some very basic notions of topology. Second, we also develop another method for analyzing and determining the nature of known critical points by resorting, when needed, to higher order derivatives. The result presented here constitutes our proposal of solution to a long standing problem: that of extending to higher dimensions what has been known for many years in the case of functions defined on open sets of the real line, regarding the classification of critical points when enough derivatives are also known at those points.*

## 1. Introduction.

The *method of Lagrange multipliers* is one of the milestones of mathematical analysis (calculus), not only because of its intrinsic, historical significance both in theoretical as well as in practical matters, but also because it has been of almost excluding use when trying to determine and analyze constrained local extrema (see, for example, [1], [2], [3], [4], [5], [11], [12], [13], [14], [15] and further references therein). In [16], F. Zizza describes two alternative methods which eliminate the multipliers from the first derivative test, i.e., from that part of the calculation in which the critical points are determined, by using differential forms. He does not indicate however a corresponding, alternative method for analyzing the critical points with regards to the possibilities of being local maxima, minima, or saddle type.

As for the analysis of critical points, one of the last articles in the literature seems to be that by D. Spring [14], where it is described how to classify constrained extrema by strongly using Lagrange multipliers. More recently, E. Constantin [4] also uses the multipliers as part of the computation, within the context of the theory of tangent sets, establishing higher order smooth necessary conditions for optimality and thus solving, in particular, some examples for which the second derivative test fails.

The first main purpose of this article is to present a comprehensive method for determining and analyzing constrained local extrema which, at the same time, represents a different alternative to the methods exposed in the first three mentioned articles, as well as in all previous works published on the topic: no use of Lagrange multipliers and/or differential forms is made in its formulation and development. Besides, it consists, basically, in reducing a constrained problem to a non-constrained one, thus fully eliminating the Lagrange multipliers from the calculations, including higher order derivatives, for all possible cases of dimension and codimension. Let us remark again here that the elimination of the Lagrange multipliers had been previously achieved, but only for the case of the first derivative test [16].

---

Partially supported by Secyt-UNC.

Moreover, it should be emphasized that the method presented here is not only very effective and straightforward, but also very simple in its conception and in the requirements of theoretical background material. In fact, in order to be performed correctly by the undergraduate students, they only need to have a basic knowledge of three very central results from calculus: the implicit function theorem, the chain rule and the Taylor's formula. Some elementary notions of topology are also needed in order to fully justify the method.

Our interest in this topic is neither limited nor new. It should be noted with respect to the latter assertion that, in fact, most of textbooks explain how to find constrained critical points by the Lagrange multipliers method, but fall short of delivering a classifying criterion similar to the one usually explained for the unconstrained case, see for example [1], [3], [5], [11], [15], amongst many others. It was precisely this absence of accurate explanation that raised the author's curiosity, since the first time he learned the subject as an undergraduate student. Thus, years later, in consultations with several fellow Mathematicians, we came to the conclusion that it was a topic worth of being elucidated, and we ourselves started by dedicating to it the corresponding parts of a book [6], now out of print, and diverse articles and congress participations afterwards. See also, for example, [7], [8], [9], [10].

Besides, the topic is not only of interest as a purely academic matter, say mainly in mathematical analysis of vector functions and teaching of mathematics at the undergraduate level, but also because of its applications to other areas, such as economics, optimization, and so on. See for instance [2], [3], [4], [12], as well as several other references therein.

With regards to the computational time needed in order to perform the operations involved, let us recall that in his article [16] F. Zizza established an experimental comparison among the Lagrange multipliers method and the two methods that he proposes as alternatives. Analogously, we make the same kind of assertions, only experimental, when it comes to deal with those matters. In fact, it should be remarked as well that the equations involved are usually nonlinear, and the appearance of Lagrange multipliers only contributes to complicate the systems even more. Thus, the question of comparing the complexity of the various algorithms involved is not as simple a matter, as only counting and comparing the number of equations in one and the others. It is a subject still due the treatment of such a complexity analysis, which goes beyond the scope of this article.

We include, as part of our exposition, carefully chosen examples where a step by step comparison is established with the first two previously mentioned methods ([16], [14]). Thus, in section 2 we consider one particular example in the lesser possible case of dimension and co-dimension, i.e., for a real valued function defined on the plane, subject to only one constraint. After some theoretical remarks, we also present alternatives to two of the examples whose solutions were previously proposed in [4].

In section 3 we develop theoretically the first of the methods proposed for the *general constraint problem*, i.e., how to determine and analyze local constrained extrema for all possible cases on the finite number of dimensions as well as co-dimensions, by proving two of the main results of this work: the first derivative test and the second derivative test, both of them without resorting to the classical Lagrange multipliers.

Two more examples, again establishing suitable comparisons with the other first two methods ([16], [14]), including a preview on the use of higher order derivatives in order to achieve the full classification in one of the examples, occupy section 4.

In section 5, we continue our exposition by remarking three important facts: in the first place we emphasize, and show, the reasons for the similarity of results between the method introduced here and F. Zizza´s second alternative method, when it comes to calculate the critical points (first derivative test). In the second place, since the recovery of Lagrange multipliers may be of importance in some applications (see for example [2], [3], [4], [12]), we show how to compute them as an additional output of the presently proposed method, but with much lesser effort and calculation time than those required when the old historical method is used. In the third place, we point out the fact that the present method does not stop with the second derivative test when it comes to analyzing critical points, as shown in two of the examples exposed in this article, but can go on, in case of failure of that test, by considering higher order derivatives until the analysis and conclusions are fully achieved, somewhat anticipating what is to follow.

In Section 6 we present the second main objective of this article: to extend to higher dimensions the use of higher order derivatives of a given function in order to classify the nature of known critical points: maxima, minima, saddle type. Recall that this is very well known for the case of dimension equal to one, i.e., for functions defined on an open set of the real numbers and, in fact, we used it in particular to independently solve previous cases: Examples 2.1, 2.2, Exercise 2.3, where the Lagrange multipliers method failed to provide the answers.

It is also to be remarked here that this is long standing problem in mathematical analysis (calculus) and that the solution proposed, requiring only that enough derivatives of a function be known at a given critical point, is achieved through the implementation of a finite sequence of constrained subsidiary problems, where just determination of, and evaluation at, critical points is required, not for the function itself but only for some of the homogeneous polynomials that appear in its Taylor´s development, and by restricting such polynomials to suitable algebraic submanifold of the corresponding unit sphere.

Most of calculations in this article were made with the help of two independent software algebraic systems: Scientific WorkPlace, Version 2.5, abbreviated from now on as Swp2.5, on one hand, and Mathematika 4.0, on the other. In each case we shall indicate which one we are using. We assume besides, without further discussion, all of the theoretical background that supports those implementations. However, a word of caution is in order here: for some real algebraic equations the systems may render solutions not belonging to the real field, but to its algebraic closure instead, the complex field with non-trivial imaginary components. Since we are only interested in the real-valued ones, those kinds of complex solutions are to be discarded from our further considerations, in every instance they may appear.

## 2. Introductory examples and theoretical considerations.

Let us consider diverse possibilities for solving the following problem:

**Example 2.1**. *Find and classify the local extrema of the function* $f(x,y)=xy$ *subject to the constraint* $G(x,y)=-2x^3+15x^2y+11y^3-24y=0$.

First, by using F. Zizza´s method [16] we have, from the expressions of the differential forms

$$df=ydx+xdy\ ,\ dG=\left(-6x^2+30xy\right)dx+\left(15x^2+33y^2-24\right)dy$$

the condition that, at the critical points, their wedge product must be vanishing, i.e.,

$$df\wedge dG=3\left(2x^3-5x^2y+11y^3-8y\right)dx\wedge dy=0\,.$$

This, together with the restriction $G(x,y)=0$, furnish the system of equations

$$\left.\begin{array}{r}2x^3-5x^2y+11y^3-8y=0\\ -2x^3+15x^2y+11y^3-24y=0\end{array}\right\},$$

whose solution set, calculated with the help of Swp2.5, is

$$\{X=(x,y)\}=\{(0,0),(1,1),(-1,-1)\}\bigcup\bar{\mathrm{S}}_4\,, \tag{2.1}$$

where $\bar{\mathrm{S}}_4$ is a set with cardinality equal to 4, i.e., having exactly four elements, described by

$$\bar{\mathrm{S}}_4=\{(x,y)\}=\left\{\left(\rho,\frac{17}{22}\rho^3-\frac{20}{11}\rho\right):\ \rho\ \ is\ a\ root\ of\ \ 17Z^4-58Z^2+32\right\}.$$

The actual values for $\rho$ are:

$$\frac{1}{17}\sqrt{\left(493+51\sqrt{33}\right)}\,,\ -\frac{1}{17}\sqrt{\left(493+51\sqrt{33}\right)}\,,\ \frac{1}{17}\sqrt{\left(493-51\sqrt{33}\right)}\,,\ -\frac{1}{17}\sqrt{\left(493-51\sqrt{33}\right)}\,.$$

Observe that, in this particular example, it would be very difficult to establish the nature of the critical points obtained. Are there among them maxima, minima, saddle type?

Second, we use the multipliers method and construct the Lagrangian function [14]:

$$L(\lambda, x, y) = f(x, y) + \lambda G(x, y) = xy + \lambda\left(-2x^3 + 15x^2 y + 11y^3 - 24y\right).$$

Then, in order to find the critical points, we have to solve the system of equations represented by $\nabla L := \left(\partial L/\partial\lambda, \partial L/\partial x, \partial L/\partial y\right) = 0$, i.e.,

$$\left.\begin{array}{r} y - 6\lambda x^2 + 30\lambda xy = 0 \\ x + 15\lambda x^2 + 33\lambda y^2 - 24\lambda = 0 \\ -2x^3 + 15x^2 y + 11y^3 - 24y = 0 \end{array}\right\}.$$

By using Swp2.5, we calculate the solution set and express it as:

$$\{(\lambda, X)\} = \left\{(0,(0,0)), \left(-\frac{1}{24},(1,1)\right), \left(\frac{1}{24},(-1,-1)\right)\right\} \bigcup \mathrm{S}_4,$$

where $\mathrm{S}_4$ is again a set with cardinality equal to 4 that can be described by

$$\mathrm{S}_4 = \{(\lambda,(x,y))\} = \left\{\left(-\frac{1}{24}\rho, \left(\rho, \frac{17}{22}\rho^3 - \frac{20}{11}\rho\right)\right) : \rho \text{ is a root of } 17Z^4 - 58Z^2 + 32\right\}.$$

This produces the same solutions as the previous method, with the additional information represented by the Lagrange multipliers at each critical point.

Next, we express the (bordered) Hessian matrix of $L$ by (see [14] )

$$HL(\lambda, X) = \begin{bmatrix} 0 & \frac{\partial G}{\partial x}(X) & \frac{\partial G}{\partial y}(X) \\ \frac{\partial G}{\partial x}(X) & \frac{\partial^2 L}{\partial x^2}(\lambda, X) & \frac{\partial^2 L}{\partial x \partial y}(\lambda, X) \\ \frac{\partial G}{\partial y}(X) & \frac{\partial^2 L}{\partial y \partial x}(\lambda, X) & \frac{\partial^2 L}{\partial y^2}(\lambda, X) \end{bmatrix} =$$

$$= \begin{bmatrix} 0 & -6x^2 + 30xy & 15x^2 + 33y^2 - 24 \\ -6x^2 + 30xy & \lambda(-12x + 30y) & 1 + 30\lambda x \\ 15x^2 + 33y^2 - 24 & 1 + 30\lambda x & 66\lambda y \end{bmatrix}.$$

To conclude the analysis, one evaluates the determinant, i.e., $\Gamma_3 = \det HL(\lambda, X)$, at the points of the solution set above and finds that the function defined by $f(x, y) = xy$, subject to the constraint $G(x, y) = -2x^3 + 15x^2 y + 11y^3 - 24y = 0$, attains local minimum values at the points of the set $\mathrm{S}_4$, and local maximum values at the points $\left(-1/24,(1,1)\right)$ *and* $\left(1/24,(-1,-1)\right)$.

However, at the critical point $(0,(0,0))$ we have $\Gamma_3 = \det HL(0,(0,0)) = 0$. So that, the Lagrange Multiplier's method does not allow any conclusion with regards to the latter point, i.e., this is an indeterminate case.

***Here is our Proposal for an Alternative Method of Solution.***
(Previously exposed in [6], [7], [8], [9], and [10]).

If we compute the Jacobian matrix of $G$,

$$G'(x,y) := \begin{pmatrix} \frac{\partial G}{\partial x} & \frac{\partial G}{\partial y} \end{pmatrix} = \begin{pmatrix} -6x^2 + 30xy & 15x^2 + 33y^2 - 24 \end{pmatrix},$$

and apply the implicit function theorem, for example at points where

$$\frac{\partial G}{\partial y} = 15x^2 + 33y^2 - 24 \neq 0,$$

we can consider the function $J(x) := f(x,y) = f(x,h(x)) = x \cdot h(x)$, where $y = h(x)$ is defined implicitly by

$$G(x,y) = -2x^3 + 15x^2 y + 11y^3 - 24y = -2x^3 + 15x^2(h(x)) + 11(h(x))^3 - 24(h(x)) \equiv 0.$$

It is easy to see that the restriction $f|_S$, *where* $S := \left\{(x,y) : G(x,y) = 0, \ \partial G/\partial y \neq 0\right\}$ has a local, constrained extremity (maximum, minimum, saddle point) at $(x_0, y_0) \in S$ if, and only if, $J$ has the same kind of (unconstrained) extremity at $x_0 \in \mathbb{R}$. Thus, we proceed to find, first, the critical points of the latter function by the usual method for the one variable case, i.e., we determine the set of points $\left\{x : J'(x) = dJ/dx = 0\right\}$. Now, in order to compute the derivative of $J$ we have, by its very definition and the chain rule, that

$$J'(x) = \frac{\partial f}{\partial x} + \frac{\partial f}{\partial y} h'(x).$$

Thus, in order to obtain the derivative of $J$ we have to compute first the derivative of $h$, and we proceed to do this implicitly: $K(x) = G(x,y) = G(x,h(x)) \equiv 0$ implies that

$$\frac{dK}{dx} = \frac{\partial G}{\partial x} + \frac{\partial G}{\partial y}\frac{dh}{dx} \equiv 0.$$

It follows from the latter that

$$h'(x) = \frac{dh}{dx} = -\frac{\frac{\partial G}{\partial x}}{\frac{\partial G}{\partial y}} = -\frac{-6x^2 + 30xy}{15x^2 + 33y^2 - 24}, \tag{2.2}$$

and we obtain for the first derivative of $J$

$$J'(x) = \frac{\partial f}{\partial x} + \frac{\partial f}{\partial y} h'(x) = y + x\left(-\frac{-6x^2 + 30xy}{15x^2 + 33y^2 - 24}\right) = \frac{-5x^2 y + 11y^3 - 8y + 2x^3}{5x^2 + 11y^2 - 8}. \tag{2.3}$$

Hence, in order to find the critical points we must solve the system of equations represented by both conditions $J'(x)=0$ and $G(x,y)=0$ which, except for the denominator in the first equation, is the same system as in the case of Zizza´s method. Swp2.5 provided us with exactly the same set of solutions, recorded above as equation (2.1).

In order to classify the critical points we compute, next, the second derivative of the function $J$, which is again related to the derivatives of $h$ up to the second order. For example, since all of the derivatives of the (constant) function $K$ must vanish, we obtain

$$\frac{d^2K}{dx^2}=\frac{\partial^2 G}{\partial x^2}+2\frac{\partial^2 G}{\partial x\partial y}\frac{dh}{dx}+\frac{\partial^2 G}{\partial y^2}\left(\frac{dh}{dx}\right)^2+\frac{\partial G}{\partial y}\frac{d^2h}{dx^2}\equiv 0 .$$

Then, it follows that

$$h''(x)=\frac{d^2h}{dx^2}=-\frac{1}{\frac{\partial G}{\partial y}}\left(\frac{\partial^2 G}{\partial x^2}+2\frac{\partial^2 G}{\partial x\partial y}h'(x)+\frac{\partial^2 G}{\partial y^2}(h'(x))^2\right)=$$

$$=-2\left(\frac{50x^5-440x^3y^2-242xy^4+352xy^2-128x-331yx^4}{(5x^2+11y^2-8)^3}\right)+$$

$$-2\left(\frac{550x^2y^3+400yx^2+605y^5-880y^3+320y}{(5x^2+11y^2-8)^3}\right)$$

and, consequently,

$$J''(x)=\frac{\partial^2 f}{\partial x^2}+2\frac{\partial^2 f}{\partial x\partial y}h'(x)+\frac{\partial^2 f}{\partial y^2}(h'(x))^2+\frac{\partial f}{\partial y}h''(x)=$$

$$=\frac{-2x}{(5x^2+11y^2-8)^3}\left(-660x^3y^2+160x^3-484xy^4+704xy^2-256x-81yx^4\right)+$$

$$+\frac{-2x}{(5x^2+11y^2-8)^3}\left(1650x^2y^3-400yx^2+1815y^5-2640y^3+960y\right)$$

Alternatively, one could obtain $J''$ by calculating directly the derivative of $J'$ from equation (2.3), and then replacing in it the value of $\frac{dy}{dx}=h'(x)$ given by equation (2.2).

Then, straightforward computations show that $J''$ is positive at all the points of the set $\bar{S}_4$. Therefore, they are local minima for the function $J$ and hence for $f$. Similarly, one verifies that, at both the points $\{(1,1)\ and\ (-1,-1)\}$, it holds $J''<0$: they are local maxima for $f$. Thus, these results coincide with the ones obtained previously by using the Lagrange multipliers method.

On the other hand, at the point $(0,0)$ it holds the condition $J''(0)=0$, so that with this method the second derivative test fails too. However, in the present situation we can consider higher-order derivatives. In fact, by following the same kind of procedure as before, when calculating the second derivative $J''$, we can go on to find the expressions of higher order derivatives. It is not difficult to compute, then, that at the origin the third derivative also vanishes, $J'''(0)=0$, but the fourth derivative is different from zero: $J^{IV}(0)=-2$. Thus, Taylor's theorem expansion of $J$ in a neighborhood of this critical point can be written $J(x) = f(x, h(x)) = -1/12\, x^4 + \cdots$, where the omitted terms involve powers in $x$ greater or equal than five. Consequently, we conclude that at the point $(0,0)$ the function $f(x,y) = xy$, subject to the constraint $G(x,y) = -2x^3 + 15x^2y + 11y^3 - 24y = 0$, reaches a local maximum.

**Remark 2.2.** At the beginning of our argument we assumed the second component of the Jacobian of $G$, $G'(x,y) := \left(\partial G/\partial x \quad \partial G/\partial y\right) = \left(-6x^2 + 30xy \quad 15x^2 + 33y^2 - 24\right)$, to be non-vanishing, i.e., $\partial G/\partial y = 15x^2 + 33y^2 - 24 \neq 0$.

Had we assumed, instead and also in order to exhaust all of possibilities, that the other component is the one different from zero, i.e., $\partial G/\partial x = -6x^2 + 30xy \neq 0$, then we could proceed in a similar fashion. Only that in this case, when applying the implicit function theorem, we would have that the condition for constraint, i.e., $G(x,y) = -2x^3 + 15x^2y + 11y^3 - 24y = 0$, now defines $x$ as a function of $y$, say $x = h(y)$. Following with the procedure we would also have that the function to optimize is now given by the expression $J(y) := f(x,y) = f(h(y), y) = h(y)\,y$, whose first derivative is

$$J'(y) = \frac{\partial f}{\partial x} h'(y) + \frac{\partial f}{\partial y} = \frac{\partial f}{\partial x}\left(-\frac{\partial G/\partial y}{\partial G/\partial x}\right) + \frac{\partial f}{\partial y} = -\frac{1}{\partial G/\partial x}\left(\frac{\partial f}{\partial x}\frac{\partial G}{\partial y} - \frac{\partial f}{\partial y}\frac{\partial G}{\partial x}\right).$$

After substituting the actual values for $f$ *and* $G$ we would find the condition

$$J'(y) = -\frac{1}{2}\frac{-5x^2y + 11y^3 - 8y + 2x^3}{-x^2 + 5xy} = 0,$$

which is totally identical with the one obtained previously, except again for the denominator, and provides, together with $G(x,y) = 0$, the same solution set as before.

Finally observe that, in order to ensure ourselves to have found all of the critical points, it is necessary to analyze both cases as above. In the general case, on the other hand, one will have to treat all possibilities where the Jacobian matrix of $G$ is non-vanishing and, in such circumstances, these may amount to some more cases to consider. Nevertheless, there will always be only a finite number to analyze and resolve.

**Remark 2.3.** We can illustrate the theoretical situation with a couple of diagrams. In Figure 1, we describe the problem we are dealing with: try to find and analyze the critical points of the real valued function in two variables $f(x,y)$, in the case where $f$ is restricted to the set

$$S := \left\{(x,y) : G(x,y) = 0,\ G'(x,y) = \left(\partial G/\partial x \quad \partial G/\partial y\right) \neq 0\right\}.$$

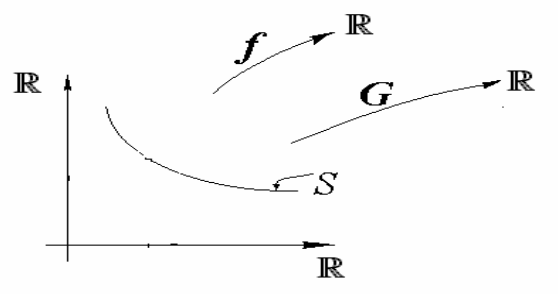

Figure 1

Let us observe, incidentally, that the components of the Jacobian matrix $G'(x,y)$ are the same as those of the gradient, $\nabla G(x,y) = \left(\partial G/\partial x, \partial G/\partial y\right)$, used by D. Spring [14], and as those of the differential form $dG = \frac{\partial G}{\partial x}dx + \frac{\partial G}{\partial y}dy$, used by F. Zizza [16] .

In Figure 2 below we describe the fact that, by assuming $\partial G/\partial y \neq 0$, the neighborhood of any given point $X_0 = (x_0, y_0) := \left(U_0, V_0\right) \in S$, with $\partial G/\partial y(X_0) \neq 0$, can be parametrized by the function $x \mapsto H(x, h(x))$ as indicated, where $y = h(x)$ is defined implicitly by the further condition $G(x,y) = G(x, h(x)) = 0$. This gives rise to the construction, by composition, of the auxiliary real-valued functions

$$J(x) := f \circ H(x) = f(x, h(x)) = f(x,y),\ K(x) := G \circ H(x) = G(x, h(x)) = G(x,y).$$

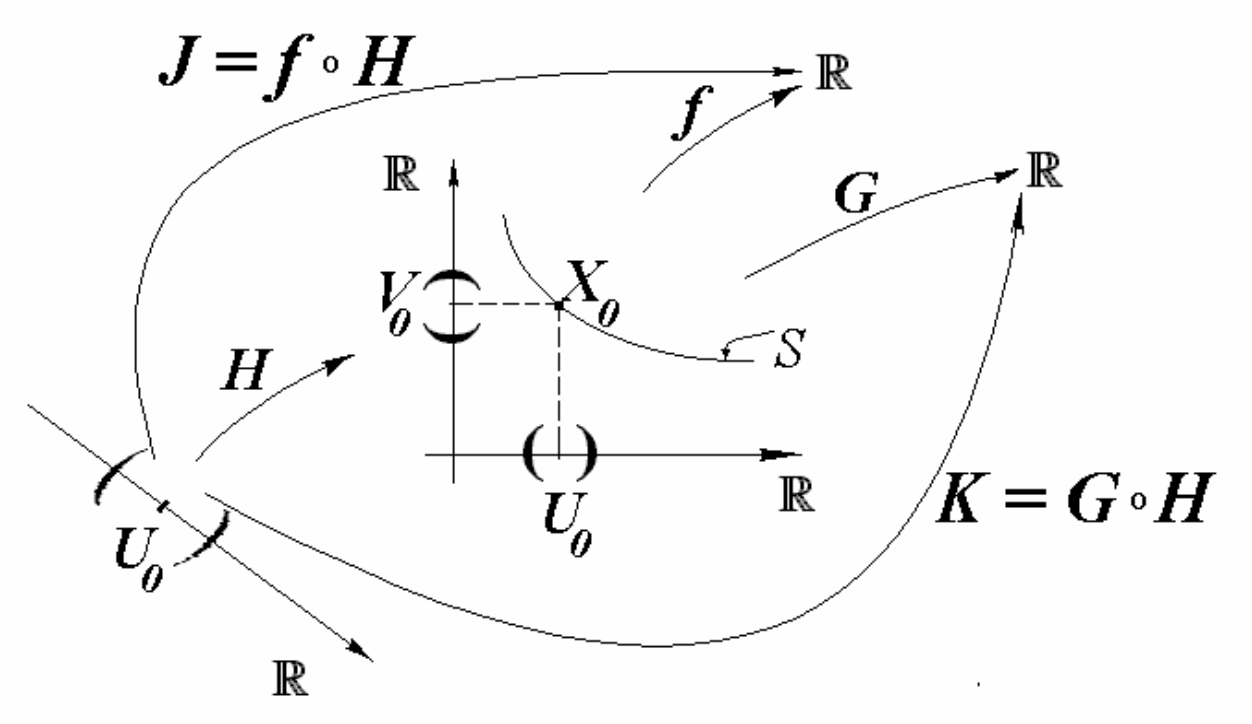

Figure 2

Henceforth, as pointed out before, the problem of finding and classifying local, constrained critical points of $f|_S$ is reduced to the same problem with respect to the (unconstrained) function $J$. However, in the expression of the latter it appears the function $h$, which may be difficult, or even impossible to calculate. Nevertheless, since we only need to obtain the derivatives of $J$, we may proceed to do so by using the chain rule, which furnishes the expression

$$J'(x) = \frac{\partial f}{\partial x} + \frac{\partial f}{\partial y} h'(x).$$

Then, in order to compute the derivative $h'(x)$, we use the fact that the function $K$ vanishes identically, so that its derivatives are also vanishing. From this we have that

$$\frac{dK}{dx} = \frac{\partial G}{\partial x} + \frac{\partial G}{\partial y} h'(x) \equiv 0$$

and, since $\partial G/\partial y \neq 0$ we obtain first, for the derivative of $h$: $h'(x) = -\frac{\partial G/\partial x}{\partial G/\partial y}$ and, for that of $J$: $J'(x) = \frac{\partial f}{\partial x} - \frac{\partial f}{\partial y}\frac{\partial G/\partial x}{\partial G/\partial y}$.

Observe that this shows, in particular, that this derivative, $J'(x)$, can always be expressed in terms of the first derivatives of the given data: $f$ and $G$.

Moreover, by the same token, if we assume enough differentiability of the data $f$, $G$, we can express the corresponding order of derivative of $J$ in terms of the derivatives of $f$ and $G$ up to the same order. For example for the second order derivative we have, since $K(x) \equiv 0$, that

$$K''(x) = \frac{d^2K}{dx^2} = \frac{\partial^2 G}{\partial x^2} + 2\frac{\partial^2 G}{\partial x \partial y}\frac{dh}{dx} + \frac{\partial^2 G}{\partial y^2}\left(\frac{dh}{dx}\right)^2 + \frac{\partial G}{\partial y}\frac{d^2 h}{dx^2} \equiv 0.$$

The latter equation implies that

$$h''(x) = \frac{d^2h}{dx^2} = -\frac{1}{\frac{\partial G}{\partial y}}\left(\frac{\partial^2 G}{\partial x^2} + 2\frac{\partial^2 G}{\partial x \partial y}\left(-\frac{\partial G/\partial x}{\partial G/\partial y}\right) + \frac{\partial^2 G}{\partial y^2}\left(-\frac{\partial G/\partial x}{\partial G/\partial y}\right)^2\right).$$

So that

$$J''(x) = \frac{\partial^2 f}{\partial x^2} - 2\frac{\partial^2 f}{\partial x \partial y}\left(\frac{\frac{\partial G}{\partial x}}{\frac{\partial G}{\partial y}}\right) + \frac{\partial^2 f}{\partial y^2}\left(\frac{\frac{\partial G}{\partial x}}{\frac{\partial G}{\partial y}}\right)^2 - \frac{\partial f}{\partial y}\left[\frac{1}{\frac{\partial G}{\partial y}}\left(\frac{\partial^2 G}{\partial x^2} - 2\frac{\partial^2 G}{\partial x \partial y}\left(\frac{\frac{\partial G}{\partial x}}{\frac{\partial G}{\partial y}}\right) + \frac{\partial^2 G}{\partial y^2}\left(\frac{\frac{\partial G}{\partial x}}{\frac{\partial G}{\partial y}}\right)^2\right)\right].$$

**Example 2.4**. Let us consider the function $F(x_1, x_2) = x_2^6 + x_1^3 + 4x_1 + 4x_2$ subject to the constraint $G(x_1, x_2) = x_1^5 + x_2^4 + x_1 + x_2 = 0$. This problem was presented as Example 2 in the article by E. Constantin [4], and solved by using the theory of tangent sets. On our hand, in the theoretical context developed above, we see that the Jacobian matrix is given by $G'(x, y) = (5x_1^4 + 1 \quad 4x_2^3 + 1)$ and we may assume first $G_2 = 4x_2^3 + 1 \neq 0$. Then, from the identity $K(x_1) = G(x_1, h(x_1)) \equiv 0$ we obtain $K'(x_1) = G_1 + G_2 h'(x_1) \equiv 0$. Thus, we have $h'(x_1) = -G_1/G_2 = -(5x_1^4 + 1)/(4x_2^3 + 1)$ and $J'(x_1) = F_1 + F_2 h'(x_1) = 3x_1^2 + 4 + (6x_2^5 + 4)(-(5x_1^4 + 1)/(4x_2^3 + 1))$. Next, we get that the system of equations represented by both conditions, $J'(x_1) = 0$ and $G(x_1, x_2) = 0$, has a solution set, obtained by using Swp2.5, with quite a few elements. In particular, the point $(0,0)$ belongs to this set. This is our primary interest here and it is easy to verify that $J'(0) = J''(0) = 0$ and $J'''(0) = 6$. It follows that $(0,0)$ is a saddle point for the problem, (*cf.* [4], pp. 47-48).

**Exercise 2.5**. In a similar fashion, consider Example 3, in the same article [4], i.e., the function $F(x_1, x_2) = x_2^6 + x_1^3 + 2x_1^2 - x_2^2 + 4x_1 + 4x_2$, subject to the same constraint as in the latter example. Show then that $J'(0) = 0$ and $J''(0) = 2$, so that $(0,0)$ is a point of strict minimum, (*cf.* [4], p. 48).

### 3. The General Constrained Problem.

In this section we shall show that extensions of the method are valid for all possible cases of larger dimension and co-dimension. In fact, in the general case we have a function

$$y = f(x_1,....,x_n,x_{n+1},....,x_{n+m}),$$

whose local extrema are to be determined in the case where the function $f$ is subject to subsidiary constraints represented by a set of equations like:

$$G_i(x_1,....,x_n,x_{n+1},....,x_{n+m}) = 0,$$

where $i = 1,2,....,m.$

We assume that the given functions are defined with enough differentiability on an open subset $U \subset \mathbb{R}^{n+m}$, $f : U \to \mathbb{R}$, $G = (G_1,...,G_m) : U \to \mathbb{R}^m$. Then, we want to determine and classify the critical points of $f$ restricted to the set $\{X = (x_1, x_2,..., x_{n+m}) : G(X) = 0\} \subset \mathbb{R}^{n+m}$. The treatment of this problem is greatly facilitated if we limit ourselves, first, to consider the case of those points where, besides, it holds the additional condition that allows applying the implicit function theorem, namely, that the rank of the corresponding Jacobian matrix is maximal, i.e., $rank(G'(X)) = m$. Thus, we consider the restriction of $f$ to the set, indeed $n$-dimensional differentiable manifold:

$$S = \{X = (x_1, x_2,..., x_{n+m}) : G(X) = 0, rank(G'(X)) = m\} \subset \mathbb{R}^{n+m}.$$

More precisely, since the function $G = (G_1,...,G_m) : U \to \mathbb{R}^m$ is assumed to be enough differentiable we can consider, first, its Jacobian matrix

$$G'(X) = \begin{bmatrix} \dfrac{\partial G_1}{\partial x_1} & \dfrac{\partial G_1}{\partial x_2} & \cdots & \dfrac{\partial G_1}{\partial x_{n+m}} \\ \vdots & \vdots & & \vdots \\ \dfrac{\partial G_m}{\partial x_1} & \dfrac{\partial G_m}{\partial x_2} & \cdots & \dfrac{\partial G_m}{\partial x_{n+m}} \end{bmatrix},$$

whose rank is assumed to equal $m$, i.e., $G'(X)$ is to have at least $m$ columns linearly independent. Let $X_0 = (x_{0_1},....,x_{0_n},x_{0_{n+1}},....,x_{0_{n+m}}) \in S$ and let us suppose for convenience, and without loss of generality, that the $m$ linearly independent columns in $G'(X_0)$ are the last ones (if not, the variables can be re-ordered): $\partial G/\partial x_{n+1}, \partial G/\partial x_{n+2},...,\partial G/\partial x_{n+m}$. Introducing the notations $U_0 = (x_{0_1},....,x_{0_n})$ and $V_0 = (x_{0_{n+1}},....,x_{0_{n+m}})$ for the components, the implicit function theorem implies that there exists a differentiable function $h = (h_1,....,h_m) : N \to \mathbb{R}^m$, defined on an open neighborhood $N \subset \mathbb{R}^n$, of $U_0$, such that it holds $h(U_0) = V_0$ and $H(u) = (u, h(u)) \in S, \forall u \in N.$ Hence, we also have that $S \cap (N \times \mathbb{R}^m) = \text{Graph}(h) = \text{Image}(H)$. In other words, the function $H : N \to \mathbb{R}^{n+m}$ realizes a vector-valued parametrization of the differentiable manifold $S$ in a neighborhood of $X_0$.

As in the case where $n = m = 1$, it is quite clear that the (constrained) function $f|_S : S \to \mathbb{R}$ reaches a local extremity (maximum, minimum or saddle point) at $X_0$ if, and only if, the (unconstrained) function $J := f \circ H : N \to \mathbb{R}$ reaches the same kind of local extremity at $U_0$. Therefore, we can reduce our original problem to the determination and analysis of the possible critical points of the latter function, following the steps required in the corresponding problem for free, unconstrained extrema.

Again, this problem would be quite easy to solve if we knew the actual, explicit expression of the function $H=(Id,h)$. However, as we also know from most of examples, this may be very difficult, or even impossible to obtain. Thus, in these cases we have to compute the derivatives of the function $h$ implicitly, by using the fact that the function $K$ obtained by composing $G$ with $H$ vanishes identically, i.e., $K:=G\circ H=G\circ(Id,h)\equiv 0$. We represent again in a graph the whole picture as Figure 3 that follows:

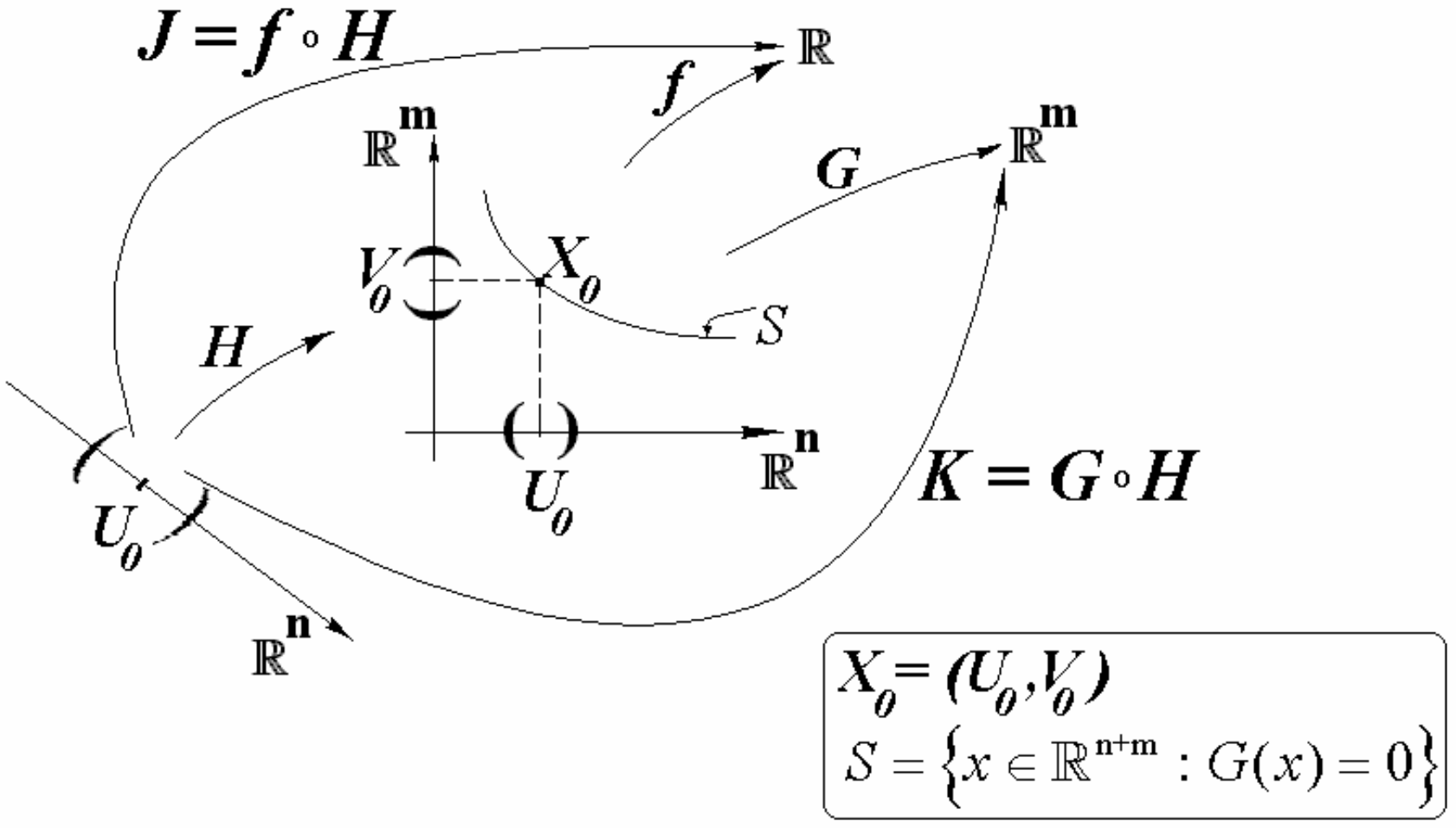

Figure 3

As stated previously we assume, for the rest of the present exposition, that both the real-valued function $f$ and the vector-valued function $G$ are enough differentiable. For example we have to require at least class $C^{(2)}$ in order to develop the second derivative test.

Both the gradient and the Hessian matrix of $J:=f\circ H$ can be computed as follows: we denote by $u=(u_1,u_2,\ldots,u_n)$ the $n$ first coordinates of a general point $X$ with coordinates $X=(x_1,....,x_n,x_{n+1},....,x_{n+m})$ and by $v=(v_1,v_2,\ldots,v_m)$ the last ones. Then, the function $h$ is defined by $h(u)=(v_1(u),v_2(u),...,v_m(u))$.

By $G_v'$ we denote the matrix whose $i-th$ row is made up with the partial derivatives of $G_i$ respect to the variables $v_1,v_2,...,v_m$, i.e.,

$$G_v'=\begin{bmatrix}\frac{\partial G_1}{\partial v_1} & \frac{\partial G_1}{\partial v_2} & \cdots & \frac{\partial G_1}{\partial v_m}\\ \vdots & \vdots & & \vdots\\ \frac{\partial G_m}{\partial v_1} & \frac{\partial G_m}{\partial v_2} & \cdots & \frac{\partial G_m}{\partial v_m}\end{bmatrix}=\begin{bmatrix}\frac{\partial G_1}{\partial x_{n+1}} & \frac{\partial G_1}{\partial x_{n+2}} & \cdots & \frac{\partial G_1}{\partial x_{n+m}}\\ \vdots & \vdots & & \vdots\\ \frac{\partial G_m}{\partial x_{n+1}} & \frac{\partial G_m}{\partial x_{n+2}} & \cdots & \frac{\partial G_m}{\partial x_{n+m}}\end{bmatrix}.$$

In what follows we denote with sub-indices the successive derivatives of the scalar and vector-valued functions needed in our exposition. Thus, for example:

$$\frac{\partial f}{\partial x_i}:=f_i;\ \frac{\partial^2 f}{\partial x_i\partial x_j}:=f_{ij};\ \frac{\partial G}{\partial x_i}:=G_i;\ \frac{\partial^2 G}{\partial x_i\partial x_j}:=G_{ij};\ \frac{\partial G_{\alpha-n}}{\partial x_i}:=G_{\alpha-n,i};\ \frac{\partial^2 G_{\alpha-n}}{\partial x_i\partial x_j}:=G_{\alpha-n,ij}.$$

Observe that, when there are already existing sub-indices which indicate components, we separate with a comma to indicate derivation of the corresponding real valued function.

Besides, we shall use the following range for indices: small Latin letters shall indicate indices running from $1$ to $n$, i.e., $1 \le i, j, k, l \ldots \le n$; while small Greek letters shall run from $n+1$ to $n+m$, i.e., $n+1 \le \alpha, \beta, \gamma, \lambda, \ldots \le n+m$. With this convention it turns out, for example, that the several expressions of the matrix $G'_v$ can indistinctly be written as

$$[G'_v] = \left( \frac{\partial G_{\alpha-n}}{\partial v_{\beta-n}} \right) = \left( \frac{\partial G_{\alpha-n}}{\partial x_\beta} \right) = \left( G_{\alpha-n,\beta} \right).$$

Moreover, since the latter is assumed to be non-singular, we shall denote its inverse by

$$[G'_v]^{-1} := \left( G^{\alpha\beta} \right).$$

Now, returning to our exposition, since the function $K$ vanishes identically on an open set, i.e., $K(u) := G(u, h(u)) \equiv 0, \forall u \in N \subset \mathbb{R}^n$, it follows that all of the derivatives of this function are also identically vanishing. In particular, for the first derivatives, respect to $u_i$, we obtain, by also using the chain rule:

$$\frac{\partial K}{\partial u_i} = \frac{\partial G}{\partial u_i}(u, h(u)) + G'_v(u, h(u)) \cdot \frac{\partial h}{\partial u_i}(u) = \begin{pmatrix} 0 \\ \cdot \\ \cdot \\ 0 \end{pmatrix},$$

which in terms of the alternative, indicial notation can also be represented by:

$$K_{\alpha-n,i} = G_{\alpha-n,i} + \sum_\beta G_{\alpha-n,\beta} h_{\beta-n,i} = 0.$$

Then,

$$\frac{\partial h}{\partial u_i}(u) = -\left[ G'_v(u, h(u)) \right]^{-1} \cdot \frac{\partial G}{\partial u_i}(u, h(u)),$$

and, again in terms of the indicial notation, we may represent by

$$h_{\gamma-n,i} = -\sum_\alpha G^{\gamma\alpha} G_{\alpha-n,i}.$$

We obtain, similarly, from the fact that the second derivatives of the function $K$ also vanish, i.e.,

$$K_{\alpha-n,ij} = G_{\alpha-n,ij} + \sum_\mu G_{\alpha-n,i\mu} h_{\mu-n,j} + \sum_\beta G_{\alpha-n,\beta j} h_{\beta-n,i} + \sum_{\beta,\nu} G_{\alpha-n,\beta\nu} h_{\beta-n,i} h_{\nu-n,j} + \sum_\beta G_{\alpha-n,\beta} h_{\beta-n,ij} \equiv 0,$$

the expression

$$h_{\gamma-n,ij} = -\sum_\alpha G^{\gamma\alpha} G_{\alpha-n,ij} + \sum_{\alpha,\mu,\rho} G^{\gamma\alpha} G_{\alpha-n,i\mu} G^{\mu\rho} G_{\rho-n,j} + \sum_{\alpha,\beta,\rho} G^{\gamma\alpha} G_{\alpha-n,\beta j} G^{\beta\rho} G_{\rho-n,i} -$$
$$- \sum_{\alpha,\beta,\nu,\rho,\mu} G^{\gamma\alpha} G_{\alpha-n,\beta\nu} G^{\nu\rho} G_{\rho-n,j} G^{\beta\mu} G_{\mu-n,i}.$$

Now, for the function $J = f \circ H$ we calculate the first derivatives $\dfrac{\partial J}{\partial u_i}(u)$:

$$J_i = f_i + \sum_{\alpha} f_\alpha h_{\alpha-n,i} = f_i - \sum_{\alpha,\beta} f_\alpha G^{\alpha\beta} G_{\beta-n,i} \,.$$

Next, we compute the derivative of this with respect to $u_j$:

$$\begin{aligned} J_{ij} = f_{ij} &- \sum_{\alpha,\beta} f_{i\alpha} G^{\alpha\beta} G_{\beta-n,j} - \sum_{\alpha,\beta} f_{\alpha j} G^{\alpha\beta} G_{\beta-n,i} + \sum_{\alpha,\beta,\mu,\nu} f_{\alpha,\beta} G^{\alpha\mu} G_{\mu-n,i} G^{\beta\nu} G_{\nu-n,j} + \\ &+ \sum_{\gamma} f_\gamma \left( \begin{aligned} &-\sum_{\alpha} G^{\gamma\alpha} G_{\alpha-n,ij} + \sum_{\alpha,\mu,\rho} G^{\gamma\alpha} G_{\alpha-n,i\mu} G^{\mu\rho} G_{\rho-n,j} + \sum_{\alpha,\beta,\rho} G^{\gamma\alpha} G_{\alpha-n,\beta j} G^{\beta\rho} G_{\rho-n,j} \\ &- \sum_{\alpha,\beta,\nu,\rho,\mu} G^{\gamma\alpha} G_{\alpha-n,\beta\nu} G^{\nu\rho} G_{\rho-n,j} G^{\beta\mu} G_{\mu-n,i} \end{aligned} \right) . \end{aligned} \tag{3.1}$$

The above theory and notations, developed so far, allow us to state the following two theorems:

**Theorem 3.1 (First Derivative Test, without multipliers).** Let $U \subset \mathbb{R}^{n+m}$ be an open set and let $f : U \to \mathbb{R}$, $G = (G_1, \ldots, G_m) : U \to \mathbb{R}^m$ be differentiable functions of class $C^{(1)}$ such that the Jacobian matrix of $G$ has maximal rank, i.e., $rank\left(G'(X)\right) = m$, and consider the restriction of $f$ to the set $S = \left\{ X = (x_1, x_2, \ldots, x_{n+m}) : G(X) = 0, rank\left(G'(X)\right) = m \right\} \subset \mathbb{R}^{n+m}$, assumed on the other hand to be non-empty and thus, indeed, an $n-$ dimensional differentiable manifold. Then, the critical points of the function $y = f(x_1, \ldots., x_n, x_{n+1}, \ldots., x_{n+m})$, subject to the constraints determined by both of conditions $G_i(x_1, \ldots., x_n, x_{n+1}, \ldots., x_{n+m}) = 0,\ i = 1, 2, \ldots., m$ and $rank\left(G'(X)\right) = m$, may be calculated, taking also into account every possible case where $m$ columns of the matrix $G'(X)$ are linearly independent, by solving, in each of those cases, the system of equations represented by

$$\left. \begin{aligned} \nabla J &= 0 \\ G &= 0 \end{aligned} \right\} .$$

**Proof.** Obvious from the above considerations, stressing once again the very important observation that, for obtaining all of the critical points one has to analyze too, separately, all of the remaining possible cases where $m$ columns of the matrix $G'(X)$ are linearly independent. Once such a listing is exhausted, the job of finding all of the critical points is done. The theorem is proved.

Observe, besides, that it may be more convenient to break the above in terms of components,

$$\left. \begin{aligned} J_i &= 0,\ i = 1, 2, \ldots, n \\ G_{\alpha-n} &= 0,\ \alpha - n = 1, 2, \ldots, m \end{aligned} \right\} .$$

Here we further observe that this is a system of $n + m$ equations in $n + m$ unknowns. In fact, in the first place, the gradient $\nabla J$ had been previously written in terms of the derivatives of the given data $f$ and $G$ up to the first order and, therefore, both the gradient $\nabla J$ and the vector-valued function $G$ can also be expressed in terms of the components of the variable, general point $X = (x_1, \ldots., x_n, x_{n+1}, \ldots., x_{n+m})$.

**Theorem 3.2 (Classical Second Derivative Test, without multipliers).** With the above notations and conditions, assume now that the functions $f$ and $G$ are differentiable of class $C^{(2)}$. Then, the critical points of the function $y = f(x_1,....,x_n,x_{n+1},....,x_{n+m})$, subject to the constraints determined by the conditions $G_i(x_1,....,x_n,x_{n+1},....,x_{n+m}) = 0,\ i = 1,2,....,m$; $rank\left(G'(X)\right) = m$ may be analyzed in order to determine whether they are maxima, minima or of saddle type, through the $n \times n$ Hessian matrix $Hess(J) := \left(J_{ij}\right)$, whose components $J_{ij}$ were described above, in equation (3.1), i.e., in terms of the derivatives of the given data $f$ and $G$ up to the second order, for each of the possible cases where $m$ columns of the matrix $G'(X)$ are linearly independent.

**Proof.** Obvious.

Finally, in case of failure of the latter, *second derivative test* (indeterminate case) we could if necessary proceed, in a similar fashion, to calculate higher order derivatives of the function $J$. Suitable use of the Taylor's formula would allow then to classify the critical points. Example 4.1, in next section, illustrates this fact, while in section 6 we present a new, comprehensive method for the analysis of critical points which, when needed, resorts to those higher order derivatives. In fact, we shall also see there that it even furnishes a very valid alternative to the latter second derivative test.

**Remark 3.3.** Insofar, when dealing with this problem, we have considered two conditions to hold: $G(X) = 0$ and the Jacobian matrix of $G$ to have maximal rank, i.e., $rank\left(G'(X)\right) = m$. However, there may exist points satisfying the first and not the second (critical points for $G$), or even it could occur that the latter condition is satisfied at no points of the (nonempty) solution set of the equation. See various examples described by J. Nunemacher in [13], and also the proof of theorem 6.3, around equation (6.6) ahead, as well as exercise 7.5. In this situation we may still try to determine and classify extrema of the given function $f$, restricted to the set $\left\{X : G\left(X\right) = 0\right\} \subset \mathbb{R}^{n+m}$, such that $rank\left(G'(X)\right) < m$. Observe, first, that the solution set of this system, assumed to be nonempty, may not contain any differentiable manifold of dimension equal to $n$; and, second, that in general the system of equations involved could be quite complicated, including even the possibility that it has transcendental terms, so that the situation may fall beyond the possibilities of treatment by means of the fields of knowledge that we have available on the topic nowadays, to obtain exact solutions. Except if, in the case of analytic functions, we consider the possibility of truncating the corresponding power series, by introducing some additional criterion for dealing with approximate solutions, question that goes far beyond the objectives and scope of the present article. However in some cases, for example if all of the equations in the system are of polynomial type and therefore the solution set is an algebraic variety, irreducible or not, we may try to resort to the possibility of expressing that solution set as a finite union of solution sets of other systems of equations such that, for every one of those new systems, the two conditions that guarantee treatment by means of the implicit function theorem are fulfilled. If that is the case, for each of those systems we should have, obviously, that the number of equations and the maximal rank, say $\overline{m}$, are less than or equal to the dimension of the ambient space, but strictly greater than the original one, i.e., $m < \overline{m} \leq m + n$, and consequently the dimension of the corresponding differential manifold equal to $m + n - \overline{m}$.

## 4. Further comparison of the three methods.

Consider, first, the following example.

**Example 4.1.** *Find and classify the local extrema of the function* $f\left(x, y, z\right) = x \cdot y \cdot z$ *subject to the constraint* $G\left(x, y, z\right) = -2x^3 + 15x^2 y + 11y^3 - 24y = 0$.

***Differential forms methods*** [16]:

From the functions $f(x,y,z)=x\cdot y\cdot z$ and $G(x,y,z)=-2x^3+15x^2y+11y^3-24y$ we obtain, by exterior differentiation:

$$df = yzdx + xzdy + xydz \ , \ dG = \left(-6x^2+30xy\right)dx + \left(15x^2+33y^2-24\right)dy + 0dz .$$

Thus, their wedge (exterior) product can be written

$$df\wedge dG = 3z\left(2x^3-5x^2y-8y+11y^3\right)dx\wedge dy + 6x^2y\left(x-5y\right)dx\wedge dz - \\ -3xy\left(5x^2+11y^2-8\right)dy\wedge dz$$

***Zizza's First alternative***.

The condition $df\wedge dG=0$, together with the constrained equation, furnish the system of four equations in three unknowns:

$$\left.\begin{array}{r} -2x^3+15x^2y+11y^3-24y=0 \\ 3z\left(2x^3-5x^2y-8y+11y^3\right)=0 \\ 6x^2y\left(x-5y\right)=0 \\ 3xy\left(5x^2+11y^2-8\right)=0 \end{array}\right\}$$

with solution set, calculated again with the help of Swp2.5, and expressed as

$$\{(0,0,z): z\in\mathbb{R}\}\bigcup\left\{\left(0,\frac{2}{11}\sqrt{66},0\right),\left(0,-\frac{2}{11}\sqrt{66},0\right)\right\}. \quad (4.1)$$

***Zizza's Second alternative***.

Consider the coordinate system $\left(y_1,y_2,y_3\right)$, where $y_1=G$, $y_2=x$ and $y_3=z$. Thus, we obtain the condition

$$dy_1\wedge dy_2\wedge dy_3 = \left(15x^2+33y^2-24\right)dx\wedge dy\wedge dz \neq 0 \text{, i.e., } 15x^2+33y^2-24\neq 0 .$$

If we now proceed to calculate the solutions of the system as suggested by Zizza, i.e.: $df\wedge dG\wedge dx=0$, $df\wedge dG\wedge dz=0$, $G=0$, we find a solution set which contains, in addition to the set determined by the previous alternative, two more points where the coordinate system fails to be valid, i.e., where $15x^2+33y^2-24=0$. Thus, in order to overcome this difficulty we consider, instead, the system represented symbolically by

$$\left.\begin{array}{c} \dfrac{df\wedge dG\wedge dx}{dy_1\wedge dy_2\wedge dy_3}=0 \\ \dfrac{df\wedge dG\wedge dz}{dy_1\wedge dy_2\wedge dy_3}=0 \\ G=0 \end{array}\right\}.$$

Observe that the quotient of differential forms makes sense since all of them are one-dimensional objects. Then, with this modified version of Zizza's Second Alternative, the solution set is the same as that obtained with the first alternative.

Finally, and again as it happened with the case of **Example 2.1**, Zizza´s method does not allow classifying any of the critical points obtained.

***Lagrange method*** **[14].**

First, we form the Lagrangian function:

$$L(\lambda,x,y,z)=f(x,y,z)+\lambda G(x,y,z)=xyz+\lambda\left(-2x^3+15x^2y+11y^3-24y\right).$$

Thus, for the present case, in order to find the critical points we have to solve the system

$$\nabla L:=\left(\partial L/\partial\lambda,\partial L/\partial x,\partial L/\partial y,\partial L/\partial z\right)=0\text{, i.e.,}$$

$$\left.\begin{aligned}-2x^3+15x^2y+11y^3-24y&=0\\ yz+\lambda\left(-6x^2+30xy\right)&=0\\ xz+\lambda\left(15x^2+33y^2-24\right)&=0\\ xy&=0\end{aligned}\right\}$$

Swp2.5 provided the solution, which we represent by

$$\{(\lambda,X)\}=\{(0,(0,0,z)):z\in\mathbb{R}\}\bigcup\left\{\left(0,\left(0,\frac{2}{11}\sqrt{66},0\right)\right),\left(0,\left(0,-\frac{2}{11}\sqrt{66},0\right)\right)\right\}.$$

Next, for the (bordered) Hessian matrix of $L$ we obtain

$$HL(\lambda,X)=\begin{bmatrix}0 & -6x^2+30xy & 15x^2+33y^2-24 & 0\\ -6x^2+30xy & \lambda(-12x+30y) & z+30\lambda x & y\\ 15x^2+33y^2-24 & z+30\lambda x & 66\lambda y & x\\ 0 & y & x & 0\end{bmatrix}.$$

Then, a straightforward computation shows that, at the critical points

$$\left(0,\left(0,\frac{2}{11}\sqrt{66},0\right)\right)\text{ and }\left(0,\left(0,-\frac{2}{11}\sqrt{66},0\right)\right)$$

we have that $\Gamma_3=0$ and $\Gamma_4=363096/1331\neq 0$. Hence, both are saddle points. However, at the rest of the critical points, i.e., $\{(0,(0,0,z)):z\in\mathbb{R}\}$, we have $\Gamma_3=0$, $\Gamma_4=0$ and, therefore, at every one of these points the method does not allow to get any conclusions, i.e., all are indeterminate cases.

***Alternative proposed method*** [6], [7], [8], [9], [10].

With the convention that the variables are relabeled as $x=u_1,z=u_2,y=v_1$, in order to suit our previously exposed notation, we consider the Jacobian matrix of the function defining the restriction $G(x,y,z)=-2x^3+15x^2y+11y^3-24y$:

$$G'(X) = \begin{bmatrix} G_1 & G_2 & G_3 \end{bmatrix} = \begin{bmatrix} -6x^2 + 30xy & 0 & 15x^2 + 33y^2 - 24 \end{bmatrix},$$

and assume that $G_3 = 15x^2 + 33y^2 - 24 \neq 0$. Hence, by the implicit function theorem, the constraint condition $G(x, y, z) = -2x^3 + 15x^2 y + 11y^3 - 24y = 0$ allows to consider $y = v_1$ as a function of $x = u_1$ and $z = u_2$, i.e., $v_1 = y = h(x, z) = h(u_1, u_2)$. Therefore, we can write

$$K(u_1, u_2) := G(u_1, u_2, h(u_1, u_2)) = -2x^3 + 15x^2 h(x, z) + 11(h(x, z))^3 - 24h(x, z) \equiv 0.$$

From the latter we obtain $K_1 = G_1 + G_3 h_1$, $K_2 = G_2 + G_3 h_2$, and it follows that

$$h_1 = -\frac{G_1}{G_3} = -2x\frac{-x + 5y}{5x^2 + 11y^2 - 8},\ h_2 = -\frac{G_2}{G_3} = \frac{0}{15x^2 + 33y^2 - 24} = 0.$$

Thus, for the first derivatives of the function $J(x, z) = f(x, h(x, z), z) = xzh(x, z)$ we find

$$J_1 = f_1 + f_2 h_1 = z\frac{-5x^2 y + 11y^3 + 2x^3 - 8y}{5x^2 + 11y^2 - 8},\ J_2 = f_2 + f_3 h_2 = xy.$$

Then, the critical points are found by solving the system

$$\left.\begin{array}{r} \dfrac{z\left(-5x^2 y + 11y^3 + 2x^3 - 8y\right)}{15x^2 + 33y^2 - 24} = 0 \\ xy = 0 \\ -2x^3 + 15x^2 y + 11y^3 - 24y = 0 \end{array}\right\}.$$

This system exhibits the same solutions as those obtained by using Zizza´s method, labeled above as (4.1). Following with the calculations we proceed to find, next, the second derivatives of the function $J$. We do this by calculating first the corresponding derivatives of $h$ which are obtained by using the fact, described previously, that all of the derivatives of $K$ vanish. In this fashion we finally obtain, by means of straightforward calculations that:

$$J_{11} = 2xz\frac{660x^3y^2 - 160x^3 + 484xy^4 - 704xy^2 + 256x}{\left(5x^2 + 11y^2 - 8\right)^3} +$$
$$+\frac{81yx^4 - 1650x^2y^3 + 400x^2y - 1815y^5 + 2640y^3 - 960y}{\left(5x^2 + 11y^2 - 8\right)^3},$$
$$J_{12} = J_{21} = \frac{-5x^2y + 11y^3 + 2x^3 - 8y}{5x^2 + 11y^2 - 8},\ J_{22} = 0.$$

It follows that, at the critical point $\left(0, \tfrac{2}{11}\sqrt{66}, 0\right)$, the Hessian matrix of $J$ is given by

$$Hess(J)=\begin{bmatrix} 0 & \frac{2}{11}\sqrt{66} \\ \frac{2}{11}\sqrt{66} & 0 \end{bmatrix},$$

and since the eigenvalues of the latter are $2/11\sqrt{66}$ and $-2/11\sqrt{66}$ there is saddle at that point, conclusion that coincides with the one obtained by using Lagrange multipliers method. Similar situation happens at the symmetric, critical point $\left(0,-2/11\sqrt{66},0\right)$.

Now, at the rest of critical points it is easy to see that the Hessian matrix vanishes, so that the second derivative test fails with this method too. But, we can follow with the computation and analysis of higher-order derivatives. In fact, by a similar procedure we find first that, at the critical points of the form $\{(0,0,z): z\in\mathbb{R}\}$, all of the third order derivatives vanish, i.e., $J_{111}=J_{222}=J_{122}=J_{112}=0$, while most of the fourth order ones are also vanishing, except one of them: $J_{1111}=-2z$. It follows, by also using Taylor's theorem expansion for $J$ in a neighborhood of each critical point, that all of the points of the form $\left(0,0,z_0\right)$ with $z_0>0$ are maxima, while those with $z_0<0$ are minima. It also follows that at the origin, $(0,0,0)$, all of the fourth order derivatives vanish, so that we have to resort to the fifth order ones. In fact, most of these are also vanishing, with the exception of those having four indices equal to 1 and one equal to 2, i.e., $J_{11112}=J_{11121}=J_{11211}=J_{12111}=J_{21111}=-2$. We conclude that the function has a saddle point at the origin. This completes the analysis of critical points by the present method. The full theoretical justification of what we have exposed here in this paragraph will take place in section 6 ahead, where the so-called **Higher Order Derivative Test** is presented.

To finish this section, let us consider the following problem, which represents an extension of the one exposed by F. Zizza in [16]:

**Example 4.2.** *Find and classify the local extrema of the function defined by*

$$f(u,v,x,y)=(x-u)^2+(y-v)^2,$$

*subject to the constraints represented by the equations*

$$G_1(u,v,x,y)=\frac{x^2}{4}+\frac{y^2}{9}-1=0,\quad G_1(u,v,x,y)=(u-3)^2+(v+5)^2-1=0\ .$$

Recall that Zizza´s problem was only to find the minimum distance between points on the ellipse and points on the circle. We seek to perform the full analysis of the problem by using the present approach.

By enumerating the variables as: $u_1=u,\ u_2=x,\ v_1=h_1=v,\ v_2=h_2=y$, we compute the Jacobian matrix of $G=\left(G_1,G_2\right)$:

$$G'(X)=\begin{bmatrix} G_{1,1} & G_{1,2} & G_{1,3} & G_{1,4} \\ G_{2,1} & G_{2,2} & G_{2,3} & G_{2,4} \end{bmatrix}=\begin{bmatrix} 0 & x/2 & 0 & 2y/9 \\ 2u-6 & 0 & 2v+10 & 0 \end{bmatrix},$$

and assume that the minor

$$G'_v=\begin{bmatrix} G_{1,3} & G_{1,4} \\ G_{2,3} & G_{2,4} \end{bmatrix}=\begin{bmatrix} 0 & 2y/9 \\ 2v+10 & 0 \end{bmatrix}$$

is non-singular. Its inverse matrix is then given by

$$[G'_v]^{-1} = \begin{bmatrix} G^{11} & G^{12} \\ G^{12} & G^{22} \end{bmatrix} = \begin{bmatrix} 0 & 1/2(v+5) \\ 9/2y & 0 \end{bmatrix}.$$

The first derivatives of $f$ are $f_1 = -2(x-u), f_2 = 2(x-u), f_3 = -2(y-v), f_4 = 2(y-v)$; from the condition $K(u_1,u_2) := G(u_1,u_2,h_1(u_1,u_2),h_2(u_1,u_2)) \equiv 0$ we obtain

$$h_{1,1} = -(u-3)/(v+5),\ \ h_{2,1} = 0,\ \ h_{1,2} = 0,\ \ h_{2,2} = -9x/4y\,;$$

and, consequently, the first derivatives of $J := f \circ H$ may be represented by

$$J_1 = f_1 + f_3 h_{1,1} + f_4 h_{2,1} = -2x + 2u - 2(-y+v)\frac{u-3}{v+5},$$
$$J_2 = f_2 + f_3 h_{1,2} + f_4 h_{2,2} = 2x - 2u + \frac{9}{2}(-y+v)\frac{x}{y}.$$

In a similar fashion we compute the second derivatives:

$$J_{11} = 2\frac{5v^2 + 50v + 170 + 5u^2 - 30u + yv^2 + 10yv + 34y + yu^2 - 6yu}{(v+5)^3},$$
$$J_{12} = J_{21} = -\frac{1}{2}\frac{4yv + 20y + 9xu - 27x}{y(v+5)},$$
$$J_{22} = \frac{1}{8}\frac{-20y^3 + 36y^2 v + 81x^2 v}{y^3}.$$

Now, the critical points are obtained by solving the system

$$\left.\begin{aligned} -2x + 2u - 2(-y+v)\frac{u-3}{v+5} &= 0 \\ 2x - 2u + \frac{9}{2}(-y+v)\frac{x}{y} &= 0 \\ \frac{x^2}{4} + \frac{y^2}{9} - 1 &= 0 \\ (u-3)^2 + (v+5)^2 - 1 &= 0 \end{aligned}\right\}$$

*Mathematika 4.0* provided the following (approximate) solution set:

$$X_1 := (u_1, x_1, v_1, y_1) = (3.41407,\ -0.580423,\ -5.91025, 2.87089),$$
$$X_2 := (u_2, x_2, v_2, y_2) = (2.58593,\ -0.580423,\ -4.08975, 2.87089),$$
$$X_3 := (u_3, x_3, v_3, y_3) = (3.64566, 0.982085,\ -5.76362,\ -2.61341),$$
$$X_4 := (u_4, x_4, v_4, y_4) = (2.35434, 0.982085,\ -4.23638,\ -2.61341).$$

This set is the same as that obtained by using F. Zizza´s first alternative. It took about the same time as that required, for solving the problem, by means of the second alternative, but without exhibiting, of course, the four points reported by that author where the coordinate system is not valid. Recall too that, in his article, he indicates the comparison in computation time with the Lagrange method.

Next, we analyze the Hessian matrix at these points. For example:

$$Hess(J)\big|_{X_1} := \begin{bmatrix} J_{11} & J_{12} \\ J_{12} & J_{22} \end{bmatrix}\Bigg|_{X_1} = \begin{bmatrix} -20.873 & -2.4139 \\ -2.4139 & -12.616 \end{bmatrix}$$

and since the (approximate) eigenvalues of this matrix are $-21.527,\ -11.962$, i.e., both are negative, it turns out that the function $J$ reaches a maximum at $X_1$. Hence, so does $f$. By a similar analysis one finds that both $X_2$ and $X_3$ are saddle points while $X_4$ is that of minimum. The (approximate) distances are: 2.1254 for the minimum; 9.647 the maximum; 7.647 and 4.1253 at the saddle points.

It is to be noted finally that, if we were to use Lagrange multiplier's method in order to classify the critical points, we would have to consider a (bordered) Hessian represented by a matrix of order $(n+2m)\times(n+2m)=6\times 6$ (see [14]).

## 5. Theoretical observations and considerations.

**Remark 5.1.** The fact that, in the examples presented here, the first derivative test of both our own approach and that of Zizza´s second alternative furnished the same critical points except for those where, in the latter, the new coordinates are not valid, is by no means a coincidence. It is not difficult to prove, in fact, that the expression for what we called the $i^{th}-partial\ derivative$ of the function $J$, i.e.,

$$J_i = f_i + \sum_{\alpha} f_\alpha h_{\alpha-n,i} = f_i - \sum_{\alpha,\beta} f_\alpha G^{\beta\alpha} G_{\beta-n,i}\,,$$

is the same as that obtained by making the quotient

$$(-1)^{i-1}\frac{df\wedge dy_1\wedge\ldots\ldots\wedge \overset{\wedge}{dy_i}\wedge\ldots\ldots dy_{n+m}}{dy_1\wedge\ldots\ldots\wedge dy_i\wedge\ldots\ldots dy_{n+m}} = (-1)^{i-1}\frac{df\wedge du_1\wedge\ldots\wedge \overset{\wedge}{du_i}\wedge\ldots\wedge du_n\wedge dG_1\wedge\ldots\wedge dG_m}{du_1\wedge\ldots\wedge du_i\wedge\ldots\wedge du_n\wedge dv_1\wedge\ldots\wedge dv_m}$$

where we have adapted and combined conveniently Zizza´s notation with ours.

**Remark 5.2.** As it was exposed, the method presented here avoids the use of Lagrange multipliers. However, in some applications, as for example in economics, optimization, etc., [2], [3], [4], [12], it may be desirable to obtain their concrete values. We show next that the historic multipliers, being no longer part of the problem, can also be obtained by the present method as an output of the solution, with no significant additional effort. In fact, suppose we have already computed the critical points, and let $X_0$ be one of them. Then, Lagrange multipliers theorem asserts that that there exist multipliers $\lambda_1,\ldots,\lambda_m$ such that $\left(f+\lambda_1 G_1+\cdots+\lambda_1 G_1\right)'(X_0)=(0,\ldots,0)$.

Since the point $X_0$ is already known, the latter represents a linear system of $n+m$ equations in the $m$ unknowns $\lambda_1,\ldots,\lambda_m$. However, by the assumptions made in the development of the method, we can reduce the above to the system of $m$ linear equations represented by

$$[\lambda_1 \ \cdots \ \lambda_m]\begin{bmatrix} \frac{\partial G_1}{\partial v_1} & \frac{\partial G_1}{\partial v_2} & \cdots & \frac{\partial G_1}{\partial v_m} \\ \vdots & \vdots & & \vdots \\ \frac{\partial G_m}{\partial v_1} & \frac{\partial G_m}{\partial v_2} & \cdots & \frac{\partial G_m}{\partial v_m} \end{bmatrix} = -\begin{bmatrix} \frac{\partial f}{\partial v_1} & \cdots & \frac{\partial f}{\partial v_m} \end{bmatrix},$$

where all derivatives are evaluated at $X_0$. Therefore, we can express the solution by

$$\lambda_{\beta-n} = -\sum_{\alpha} f_{\alpha-n}\left(X_0\right) G^{\alpha\beta}\left(X_0\right).$$

Let us apply the above to obtain the multipliers in example 4.2, which is the only remaining case where the multipliers are not known yet. We do this by specializing at each critical point the expression:

$$-[f_3 \quad f_4]\begin{bmatrix} G^{11} & G^{12} \\ G^{12} & G^{22} \end{bmatrix} = -\begin{bmatrix} -2(y-v) & 2(y-v) \end{bmatrix}\begin{bmatrix} 0 & 1/2(v+5) \\ 9/2y & 0 \end{bmatrix} =$$

$$= -\begin{bmatrix} \frac{9(y-v)}{y} & \frac{-(y-v)}{(v+5)} \end{bmatrix}.$$

Then, by substituting into the latter expression the corresponding values, already computed, we find:

At $X_1$ the multipliers are (approximately) : $\lambda_1 = -27.528, \ \lambda_2 = -9.647$

At $X_2$ : $\lambda_1 = -21.821, \ \lambda_2 = 7.647$

At $X_3$ : $\lambda_1 = 10.849, \ \lambda_2 = -4.1254$

At $X_4$ : $\lambda_1 = 5.5891, \ \lambda_2 = 2.1254$

It is to be remarked again that the whole procedure performed in this way, i.e., within the framework of the method exposed here, including the calculations of Lagrange multipliers, is, at least experimentally, much faster than trying to execute the problem and obtain its solution by the traditional Lagrange multipliers method.

**Remark 5.3.** As shown by example 4.1, the second derivative test may fail, either by using Lagrange multipliers method or the second derivative test (Theorem 3.2 above). However, as also exposed in those same examples, the use of higher order derivatives may allow to continue the further analysis of critical points, as we are about to prove in the next section.

## 6. The Higher Derivative Test.

In this section we expose another method that in case of failure of the *second derivative test,* indeterminate case, can be used in order to complete the classification of already known critical points. It consists, basically and when necessary, in calculating higher order derivatives of the function $f$. Suitable use of the Taylor's formula allows then to classify the critical points. The only necessary condition is that the latter be non trivial, i.e., that there exists at least one derivative, greater or equal than two, which is different from zero at a given critical point. Example 4.1 illustrated this quite well and the following results justify that procedure from the theoretical point of view. Observe first of all, from the above expositions and considerations, that it is enough to consider the case where the objects to analyze are the extrema of a function $f : U \to \mathbb{R}$ where $U$ is an open subset of $\mathbb{R}^d$, since we have already reduced the analysis of critical points of a constraint problem to a non-constrained one.

On the other hand, observe too that the method is very well known for the case of dimension $d=1$ and, in particular, we used it in order to solve the cases presented as Examples 2.1, 2.4 and Exercise 2.5, where the Lagrange multipliers method failed to provide an answer. Thus, we assume in this section that $d \geq 2$. It will also be clear from what follows that the method provides, besides, a valid alternative even for the classical second derivative test.

Let $X_0 = (x_{0_1}, \ldots, x_{0_d}) \in U$ be a critical point of the given function $f : U \to \mathbb{R}$ and assume that the first non-vanishing derivative at that point is of $k^{th}$-order. Then, for points $X = (x_1, \ldots, x_d)$ in a neighborhood of that point we may write, by using Taylor's formula, that

$$f(X) = f(X_0) + \frac{1}{k!}\sum f_{i_1 i_2 \ldots i_k}(X_0)\left(x_{i_1} - x_{0_{i_1}}\right)\left(x_{i_2} - x_{0_{i_2}}\right)\ldots\left(x_{i_k} - x_{0_{i_k}}\right) + \text{higher degree terms.}$$

The last displayed term is a homogeneous polynomial function of $k^{th}$-degree, represented as

$$P_k = P_k\left(\left(x_1 - x_{0_1}\right),\left(x_2 - x_{0_2}\right),\ldots,\left(x_d - x_{0_d}\right)\right) = \frac{1}{k!}\sum f_{i_1 i_2 \ldots i_k}(X_0)\left(x_{i_1} - x_{0_{i_1}}\right)\left(x_{i_2} - x_{0_{i_2}}\right)\ldots\left(x_{i_k} - x_{0_{i_k}}\right).$$

For this kind of homogeneous polynomial function it is well-known the property that, for any real number $t \neq 0$, it holds the relationship

$$P_k\left(t\left(x_1 - x_{0_1}\right),\ldots,t\left(x_d - x_{0_d}\right)\right) = t^k P_k\left(\left(x_1 - x_{0_1}\right),\left(x_2 - x_{0_2}\right),\ldots,\left(x_d - x_{0_d}\right)\right) \tag{6.1}$$

The method will rely on evaluating, on the corresponding unit sphere, the successive nonvanishing homogeneous polynomials that appear in the Taylor´s development of $f$. Thus, we consider next the:

## First Subsidiary Problem:

Determine the critical points of $P_k(X) = P_k\left(\left(x_1 - x_{0_1}\right),\left(x_2 - x_{0_2}\right),\ldots,\left(x_d - x_{0_d}\right)\right)$ subject to the restriction that $X$ belongs to the unit sphere centered at $X_0$, i.e., $X \in S^{d-1} := S^{d-1}(X_0;1) =$ $\left\{X : \left(x_1 - x_{0_1}\right)^2 + \left(x_2 - x_{0_2}\right)^2 + \cdots + \left(x_d - x_{0_d}\right)^2 = 1\right\}$.

Let us observe, first, that in the last expression, an also in what follows, we are assuming that the given sphere, border of the corresponding unit ball, is fully contained in the open set $U$, by performing a rescaling of coordinates, in case the latter were necessary in order to achieve that particularly convenient, comparative situation.

Then, if we define $G_1(X) := \left(x_1 - x_{0_1}\right)^2 + \left(x_2 - x_{0_2}\right)^2 + \cdots + \left(x_d - x_{0_d}\right)^2 - 1$ the restriction condition is given by $G_1(X) = 0$, with Jacobian $G_1'(X) = 2\left(x_1 - x_{0_1}, x_2 - x_{0_2}, \ldots, x_d - x_{0_d}\right)$. It follows, by considering all of the possibilities where one of these components is non-vanishing, i.e., $x_i - x_{0_i} \neq 0, i = 1, \ldots, d$, that the sphere $S^{d-1}$ is covered by a finite number of the corresponding parametrizing functions, all of them in the form of Monge, i.e., graph maps, as indicated by the method exposed in the proof of Theorem 3.1. Indeed, there are $d$ of them. However, if one requires connectedness for the image of those functions, that number rises to $2d$, $2$ for each of the coordinate maps. As it is very well-known, this is one of the most common and useful ways of parametrization for the sphere, and we shall refer to it as the *distinguished or canonical parametrization.*

In other words, the initial problem to be considered is that we have a homogeneous polynomial function $P_k : \mathbb{R}^d \to \mathbb{R}$ and want to find all of the critical points of $P_{k|S^{d-1}} : S^{d-1} \to \mathbb{R}$. This is a continuous function and, hence, there must be at least two points of extreme on the (compact) set $S^{d-1}$: an absolute maximum $X_{Max}$ and an absolute minimum $X_{\min}$, which could reduce to a single point, $X_{Max} = X_{\min}$, in case both coincide. These considerations guarantee that some of the systems of equations below, which were indicated in the formulation of the method in order to determine the critical points, exposed in section 3, have real solutions

$$\left.\begin{array}{r}\nabla J_{P_k}(X)=0\\ G_1(X)=0\end{array}\right\}, \tag{6.2}$$

with $J_{P_k} = P_k \circ H$. Observe besides that, whereas we have stressed the dependence of such a function on $P_k$, it should also be remarked once again that it depends on the parametrization chosen as well.

Now, since the latter is a polynomial system of equations there could exist, in the best of possible situations, only a finite number of solutions $X_1, X_2, \ldots, X_{k_0} \in S^{d-1}$ because, as pointed out before, there should exist at least two: an absolute maximum $X_{Max}$ and an absolute minimum $X_{\min}$. However, the solution set could also consists of a finite union of objects with diverse dimensions, some of them being non-finite subsets of the sphere $S^{d-1}$ which, in a such a case, will turn out to be, as we shall see ahead, algebraic differentiable submanifolds of $S^{d-1}$ with dimensions strictly less than $\dim\left(S^{d-1}\right) = d-1$.

Next, by computing the values of $P_k$ at the solution set, for example in case there are only a finite listing of points, $P_k\left(X_1\right), P_k\left(X_2\right), \ldots, P_k\left(X_{k_0}\right)$, we find by comparison that

$$P_k\left(X_{\min}\right) \le P_k\left(X_j\right) \le P_k\left(X_{MAX}\right), j = 1, 2, \ldots, k_0.$$

Moreover, we also have that $P_k\left(X_{\min}\right) \le P_k\left(X\right) \le P_k\left(X_{MAX}\right)$, for every $X \in S^{d-1}$.

**Remark 6.1.** Let us observe that the image set $P_k\left(S^{d-1}\right)$ could reduce to a single point. This happens, for example, in the case where

$$P_k\left(X\right) = P_{2p}\left(X\right) = \lambda\left(\left(x_1 - x_{0_1}\right)^2 + \left(x_2 - x_{0_2}\right)^2 + \cdots + \left(x_d - x_{0_d}\right)^2\right)^p, \tag{6.3}$$

the image being precisely the real number $\lambda$, which must be non-vanishing since we are assuming $P_{2p} = P_k \neq 0$. It will also be useful, and very convenient in our future exposition, to consider a kind of converse to this result, in the form of various equivalent properties. Furthermore, in order to make things easier for such a purpose, let us assume, without loss of generality, that the center lies at the origin of coordinates, i.e., $X_0 = (x_{0_1}, x_{0_2}, \ldots, x_{0_d}) = (0, 0, \ldots, 0) = 0$ by performing if necessary a translation of coordinates:

**Lemma 6.2.** Let $P_k : \mathbb{R}^d \to \mathbb{R}$, with $P_k = P_k\left(x_1, x_2, \ldots, x_n\right) = \sum a_{i_1 i_2 \ldots i_k} x_{i_1} x_{i_2} \ldots x_{i_k}$, be a non-trivial homogeneous polynomial function. Then, the following statements are equivalent:

1) The homogeneous polynomial function $P_k$ restricted to $S^{d-1}$ is identically equal to a non-vanishing constant, $P_{k|_{S^{d-1}}} \equiv \lambda \neq 0,\ \lambda \in \mathbb{R}$.

2) For every one of the distinguished parametrizations of the sphere $S^{d-1}$, defined above, the

corresponding composite function $J_{P_k} = P_k \circ H$ is also identically equal to the same non-vanishing constant, $J_{P_k} = P_k \circ H \equiv \lambda$.

3) All of the first order derivatives of the latter functions are identically vanishing, in their respective domains of definition.

4) The polynomial function $P_k : \mathbb{R}^d \to \mathbb{R}$ is the one represented as equation (6.3) above, with $X_0 = 0$. i.e., there exists a positive integer $p \in \mathbb{Z}_{>0}$ and a non-vanishing real number $\lambda \in \mathbb{R},\ \lambda \neq 0$ such that $P_k(X) = P_{2p}(X) = \lambda\left(x_1^2 + x_2^2 + \cdots + x_d^2\right)^p$.

**Proof**. Equivalence of the first three statements is rather obvious.

Let us prove next that 4) implies 3). In this context, the Jacobian becomes now $G_1'(X) = 2(x_1, x_2, ..., x_d)$ and if we consider the parametrization determined by the case where a given $x_i \neq 0$ we obtain, for every $j \neq i$, that the derivative with respect to $x_j$ vanishes identically, i.e.,

$$J_{P_k,j} = \frac{\partial J_{P_k}}{\partial x_j} = \frac{\partial P_k \circ H}{\partial x_j} =$$

$$= 2\lambda p\left(x_1^2 + \cdots + x_d^2\right)^{p-1} x_j - 2\lambda p\left(x_1^2 + \cdots + x_d^2\right)^{p-1} x_i \left(\frac{x_j}{x_i}\right) \equiv 0.$$

In order to finish the proof, let us prove that 1) implies 4). And let us observe first that the degree $k$ cannot be odd, because of equation (6.1). In fact, if we assume it to be odd, by putting $t = -1$ in that equation, it leads to a contradiction. Thus $k$ is even, i.e., there exists a positive integer $p \in \mathbb{Z}_{>0}$ such that $k = 2p$.

Then, we may proceed by induction on the dimension $d$. Thus, let us take first, for $d = 2$, a homogeneous polynomial function $P_k : \mathbb{R}^2 \to \mathbb{R}$ such that its restriction to $S^1$ is identically equal to a nonvanishing constant, $P_{k|S^1} \equiv \lambda \neq 0,\ \lambda \in \mathbb{R}$. We may display the homogeneous polynomial function as follows:

$$P_k(X) = P_{2p}(X) = a_1^k x_1^{2p} + a_{12}^{2p-1,1} x_1^{2p-1} x_2 + a_{12}^{2p-2,2} x_1^{2p-2} x_2^2 + a_{12}^{2p-3,3} x_1^{2p-3} x_2^3 + a_{12}^{2p-4,4} x_1^{2p-4} x_2^4 + $$
$$+ \cdots + a_{12}^{2p-i,i} x_1^{2p-i} x_2^i + \cdots + a_{12}^{1,2p-1} x_1 x_2^{2p-1} + a_2^k x_2^{2p},$$

with $X = (x_1, x_2)$ and where we have abbreviated conveniently the notation for the constant coefficients: the subscript indices running in this case only from $1$ to $2$, while the superscripts indicate the number of times that those are repeated, separating by a comma when necessary. Thus, for example, $a_{12}^{3,4} := a_{1112222}$, $a_{12}^{4,0} = a_1^4 := a_{1111}$, $a_{12}^{0,5} = a_2^5 := a_{22222}, \ldots$

We may choose to work within the parametrization defined by assuming $x_2 \neq 0$ and compute, next, the first derivative of the composite function $J_{P_k} = P_k \circ H$:

$$J_{P_k,1} = 2pa_1^k x_1^{2p-1} + (2p-1)a_{12}^{2p-1,1} x_1^{2p-2} x_2 + (2p-2)a_{12}^{2p-2,2} x_1^{2p-3} x_2^2 +$$
$$+ \cdots + (2p-i)a_{12}^{2p-i,i} x_1^{2p-i-1} x_2^i + \cdots + a_{12}^{1,2p-1} x_2^{2p-1} +$$
$$+ \left(a_{12}^{2p-1,1} x_1^{2p-1} + 2a_{12}^{2p-2,2} x_1^{2p-2} x_2 + \cdots + i a_{12}^{2p-i,i} x_1^{2p-i} x_2^{i-1} + \cdots + 2pa_2^k x_2^{2p-1}\right)\left(-\frac{x_1}{x_2}\right).$$

By collecting terms we find:

$$J_{P_k,1} = \left(2pa_1^k - 2a_{12}^{2p-2,2}\right)x_1^{2p-1} + \left((2p-1)a_{12}^{2p-1,1} - 3a_{12}^{2p-3,3}\right)x_1^{2p-2}x_2 + \\ + \cdots + \left((2p-i)a_{12}^{2p-i,i} - (i+2)a_{12}^{2p-i-2,i+2}\right)x_1^{2p-i-1}x_2^i + \cdots + \\ + \left(a_{12}^{2,2p-2} - 2pa_2^k\right)x_1x_2^{2p-2} + a_{12}^{1,2p-1}x_2^{2p-1} - a_{12}^{2p-1,1}\frac{x_1^{2p}}{x_2}.$$

Thus, the condition $J_{P_k,1} \equiv 0$ implies successively, first that

$$a_{12}^{1,2p-1} = a_{12}^{2p-1,1} = 0 = a_{12}^{2p-3,3} = a_{12}^{2p-5,5} = \cdots = a_{12}^{2p-(2l+1),2l+1} = a_{12}^{3,2p-3} = \\ = a_{12}^{5,2p-5} = \cdots = a_{12}^{2l+1,2p-(2l+1)} = 0.$$

And then, that

$$2pa_1^k - 2a_{12}^{2p-2,2} = 0 \Rightarrow a_{12}^{2p-2,2} = pa_1^k.$$

$$(2p-2)a_{12}^{2p-2,2} - 4a_{12}^{2p-4,4} = 0 \Rightarrow a_{12}^{2p-4,4} = \frac{2p-2}{4}a_{12}^{2p-2,2} = \frac{p(p-1)}{2}a_1^k.$$

$$(2p-4)a_{12}^{2p-4,4} - 6a_{12}^{2p-6,6} = 0 \Rightarrow a_{12}^{2p-6,6} = \frac{p-2}{3}a_{12}^{2p-4,4} = \frac{p(p-1)(p-2)}{3!}a_1^k.$$

$$(2p-6)a_{12}^{2p-6,6} - 8a_{12}^{2p-8,8} = 0 \Rightarrow a_{12}^{2p-8,8} = \frac{p-3}{4}a_{12}^{2p-6,6} = \frac{p(p-1)(p-2)(p-3)}{4!}a_1^k = \frac{p!}{4!(p-4)!}a_1^k.$$

. . . . . . . . . . . . . . . . . . . . . . . . . . . . . . . . . . . . . . . . . . . . . . . . . . . . . . . . . . . . . . . . . . . . . .

$$(2p-2l)a_{12}^{2p-2l,2l} - (2l+2)a_{12}^{2p-2l-2,2l+2} = 0 \Rightarrow a_{12}^{2p-2l-2,2l+2} = \frac{p-l}{l+1}a_{12}^{2p-2l,2l} = \frac{p!}{(l+1)!(p-l-1)!}a_1^k.$$

. . . . . . . . . . . . . . . . . . . . . . . . . . . . . . . . . . . . . . . . . . . . . . . . . . . . . . . . . . . . . . . . . . . . . . . . . . . . . .

Thus, for $l = p-1$ we obtain

$$a_2^k = a_{12}^{0,2p} = \frac{1}{p}a_{12}^{2,2p-2} = \frac{p!}{p!0!}a_1^k = a_1^k.$$

Since, obviously, we also have that

$$\lambda = P_k(1,0) = P_k(0,1) = a_2^k = a_1^k,$$

it follows that

$$P_k(x_1,x_2) = P_{2p}(x_1,x_2) = \lambda\left(x_1^2 + x_2^2\right)^p,$$

and the result is valid for $d = 2$.

To finish the proof, let us assume that that the result holds in every case where the dimension of the sphere equals $d-2$ and assume 1), i.e., that $P_k = P_k\left(x_1, x_2, ..., x_d\right) = \sum a_{i_1 i_2 ... i_k} x_{i_1} x_{i_2} ... x_{i_k}$ is equal to a non-vanishing constant $\lambda \neq 0,\ \lambda \in \mathbb{R}$ when restricted to $S^{d-1}$. Observe, first of all, that

$$\lambda = P_k\left(e_i\right) = P_{2p}\left(e_i\right) = a_i^{2p} = a_i^k,\ \text{for every } i = 1, 2, ..., d\ ,$$

where $e_i = \left(0, 0, ..., 0, 1, 0, ..., 0\right)$ represent the elements of the usual Euclidean basis of $\mathbb{R}^d$.

Next, we consider the further restriction of $P_{k|S^{d-1}}$ to the $\left(d-2\right)$ – dimensional sphere obtained by intersecting $S^{d-1}$ with the hyperplane $\Pi_d := \left\{X = \left(x_1, x_2, ..., x_d\right) : x_d = 0\right\}$, $S^{d-2} = S^{d-1} \cap \Pi_d$. Obviously, we still have $P_{k|S^{d-2}} \equiv \lambda \neq 0,\ \lambda \in \mathbb{R}$, so that we may write, by inductive hypothesis, that

$$P_{k|\Pi_d} = P_k\left(x_1, x_2, ..., x_{d-1}, 0\right) = \lambda\left(x_1^2 + x_2^2 + \cdots + x_{d-1}^2\right)^p .$$

It follows that we can also write

$$P_k = P_k\left(x_1, x_2, ..., x_{d-1}, x_d\right) = \lambda\left(x_1^2 + x_2^2 + \cdots + x_{d-1}^2\right)^p + \sum_{i_j = d} a_{i_1 i_2 ... i_k} x_{i_1} x_{i_2} ... x_{i_k} ,$$

where $\sum_{i_j = d} a_{i_1 i_2 ... i_k} x_{i_1} x_{i_2} ... x_{i_k}$ is a homogeneous polynomial function of degree equal to $k$ with each term containing at least one of the variables equal to the last one in the listing, $x_{i_j} = x_d$. In order to determine the coefficients in the latter polynomial we consider the restriction of $P_k$ to the plane spanned by $e_d$ and any of the other members of the usual basis, say $e_i$, $i \neq d$. By denoting this plane by $\Pi_{i,d}$, we still have that $P_{k|S^{d-1} \cap \Pi_{i,d}} \equiv \lambda \neq 0,\ \lambda \in \mathbb{R}$ and, by taking into account the proof for the case of dimension $d = 2$ exposed above, we can also write for the restriction to the one dimensional sphere $S_{i,d}^1 := S^{d-1} \cap \Pi_{i,d}$ that

$$P_{k|S_{i,d}^1} = P_k\left(0, ..., 0, x_i, 0, ..., 0, x_d\right) = \lambda\left(x_i^2 + x_d^2\right)^p .$$

Hence, we finally obtain

$$P_k\left(x_1, x_2, ..., x_d\right) = \lambda\left(\left(\sum_{i=1}^{d-1} x_i^2\right)^p + p\left(\sum_{i=1}^{d-1} x_i^2\right)^{p-1} x_d^2 + \cdots + \frac{p!}{(p-l)!\,l!}\left(\sum_{i=1}^{d-1} x_i^2\right)^{p-l} x_d^{2l} + \cdots + x_d^{2p}\right) =$$

$$= \lambda\left(\left(\sum_{i=1}^{d-1} x_i^2\right) + x_d^2\right)^p = \lambda\left(x_1^2 + x_2^2 + \cdots + x_d^2\right)^p ,$$

and the lemma is proved.

As we shall see, it is very easy to obtain conclusions when the equivalent conditions of the latter lemma hold. Thus, including this particular possibility as a sub-case, we state the next auxiliary result.

**Lemma 6.3.** Let $P_k : \mathbb{R}^d \to \mathbb{R}$ be a nontrivial homogeneous polynomial function. Then, it follows that:

a) The image set $P_k\left(S^{d-1}\right)$ is a closed interval in the real line, i.e., $P_k\left(S^{d-1}\right)=[a,b]\subset\mathbb{R}$, with $a=P_k\left(X_{\min}\right)\le b=P_k\left(X_{MAX}\right)$.

b) If $k$ is odd, then $a<0<b$.

c) In case $k$ is even we have the following possible situations regarding the extreme values of the interval:

c$_1$) $0<a\le b$;

c$_2$) $a\le b<0$;

c$_3$) $a<0<b$;

c$_4$) $0=a<b$;

c$_5$) $a<b=0$.

**Remark 6.4.** it is easy to construct examples of those five possible cases. So this is left to the reader as an exercise.

**Proof**. We treat separately the various cases:

a) The sphere $S^{d-1}$ is a compact and connected subset of $\mathbb{R}^d$. Then its image under the (continuous) function $P_k$ must also be compact and connected as a subset of $\mathbb{R}$, thus a closed interval $[a,b]$.

b) Since we assumed that $P_k$ is nontrivial, then there exists a point $X\in S^{d-1}$ such that $P_k(X)\neq 0$. It follows from (6.1) that $P_k(-X)=(-1)^k P_k(X)=-P_k(X)$, which proves our assertion.

c) It is obvious that the stated ones cover all of the different possibilities for the closed interval $[a,b]$, with $a\le b$.

<u>**Theorem 6.5. (Higher Derivative Test)**</u>. Let $X_0=(x_{0_1},....,x_{0_d})\in U$ be a critical point of the function $f:U\to\mathbb{R}$, $U$ open subset of $\mathbb{R}^d$, and assume that the first nonvanishing derivative at that point is of $k^{th}$-order, with $k\ge 2$. Then

1) $k$ odd implies that $X_0$ is a saddle point.

2) For $k$ even we have, for the first non-vanishing, homogeneous polynomial function $P_k$ the five possibilities described in the previous lemma and, according to them, it follows that:

In case c$_1$), $f\left(X_0\right)$ is a strict local minimum.

In case c$_2$), $f\left(X_0\right)$ is a strict local maximum.

In case c$_3$), $X_0$ is a saddle point.

In case c$_4$), $f\left(X_0\right)$ is a candidate to furnish a minimum. However, there exists a subset of $S^{d-1}$, $\left(P_k\right)^{-1}(0)\cap S^{d-1}\subset S^{d-1}$, symmetric with respect to the origin, with not necessarily a finite number of points, where the $k^{th}$-homogeneous polynomial $P_k$ vanishes and we need a further, deeper analysis. In fact, the algebraic set $\left(P_k\right)^{-1}(0)\cap S^{d-1}$ may be expressed as a finite union of algebraic submanifolds of the sphere $S^{d-1}$, with dimensions strictly less than $d-1$, and one may need to repeat the analysis on the behavior of the next nonvanishing homogeneous polynomial, say $P_l, l>k$, on each one of those submanifolds. By considering the expression of the Taylor's formula in every one of the directions determined by that set, we have the following sub-cases that we label as c$_{41}$), c$_{42}$), c$_{43}$), c$_{44}$), c$_{45}$), and exhibit next:

$c_{41})$ If in the Taylor´s development of $f$ around $X_0$ there are no further homogeneous polynomials, or if there are further homogeneous polynomials but are all vanishing on the set $(P_k)^{-1}(0) \cap S^{d-1}$, it follows that $X_0$ is a point of non-strict local minimum for $f$. Besides, if there are only a finite number of homogeneous polynomials that vanish on the latter set, but there exists a next following one non-vanishing, one proceeds to apply to that polynomial the considerations that follow. In such a case, in order to avoid unnecessary complications in notation, we still call $P_l$ the corresponding homogeneous polynomial.

$c_{42})$ If the next nonvanishing homogeneous polynomial, say $P_l, l > k$, is of odd order and there exists a point $X_1 \in (P_k)^{-1}(0) \cap S^{d-1}$ such that $P_l(X_1) \neq 0$, it follows that $X_0$ is a saddle point.

$c_{43})$ If the next nonvanishing homogeneous polynomial $P_l$ is of even order and $P_l(X) > 0$, for every $X \in (P_k)^{-1}(0) \cap S^{d-1}$, it follow that $X_0$ is a point of strict local minimum.

$c_{44})$ If the next nonvanishing homogeneous polynomial $P_l$ is of even order, but there exists some point $X_1 \in (P_k)^{-1}(0) \cap S^{d-1}$ such that $P_l(X_1) < 0$, then $X_0$ is a saddle point.

$c_{45})$ If the next nonvanishing homogeneous polynomial $P_l$ is of even order and $P_l(X) \geq 0$, for every $X \in (P_k)^{-1}(0) \cap S^{d-1}$, with strict inequality at some points, but with the existence of at least one point $X_1 \in (P_k)^{-1}(0) \cap S^{d-1}$ such that $P_l(X_1) = 0$, then it is kept the expectation that $X_0$ be a local minimum. However, in that case one needs to reiterate the latter analysis resorting, if necessary, to the corresponding further subsidiary problems, by considering first the next nonvanishing homogeneous polynomial, beyond $P_l$, appearing in the Taylor´s development of the function $f$ at the critical point $X_0$ restricted, in this occasion, to the algebraic set $(P_l)^{-1}(0) \cap (P_k)^{-1}(0) \cap S^{d-1}$. Again, that set may be expressed as a finite union of algebraic submanifolds of the sphere $S^{d-1}$, with dimensions strictly less than $d-2$.

In case $c_5)$ the corresponding reverse situation to that in case $c_4)$ holds, i.e., $X_0$ may be a saddle point; or a point of strict local maximum; or a point of non-strict local maximum.

Finally, since in most of cases dimensions strictly diminish with each new nonvanishing homogeneous polynomial analyzed, except possibly in cases like $c_{41})$, we conclude that the procedure gets exhausted after a finite number of steps.

**Proof**. We consider one by one the cases stated:

1) By the proof of b) in the previous lemma, we have a point $X \in S^{d-1}$ such that $P_k(X) \neq 0$, and we consider the values of $f$ along the straight line $\{X_0 + t(X - X_0) : t \in \mathbb{R}\}$. We obviously have that

$$f(X_0 + t(X - X_0)) = f(X_0) + \alpha t^k + \text{higher degree terms in } t,\ \alpha \neq 0.$$

Since we are assuming that $k$ is odd, it follows that at every neighborhood of $X_0$ there are points, say $Y, Z$, such that $f(Y) < f(X_0) < f(Z)$, proving our assertion.

2) Next, we assume that $k$ is even and consider first the case described as $c_1)\ 0 < a \leq b$. Then it follows that, for every $X \in S^{d-1}$, $P_k(X) > 0$. By analyzing the values of the function $f$, as in the previous case, we can assert that there exists a positive real number, $r > 0$, such that for every $t$ with $|t|<r$ it holds $f\left(X_0+t\left(X-X_0\right)\right)>f\left(X_0\right)$. By continuity of $f$ we conclude that there exists another positive real number, $\delta_r>0$, such that for every point $Y$ within the open cylinder with axis $\left\{X_0+t\left(X-X_0\right): t\in(-r,r)\right\}$ and radius $\delta_r$ it holds $f(Y)>f\left(X_0\right)$. By projecting the cylinder along its axis onto the sphere $S^{d-1}$ we obtain an open $d-1$ (geodesic) ball around $X$, $B_X$. The union of all of those balls, $\bigcup_{X\in S^{d-1}} B_X$, provides an open covering for the sphere $S^{d-1}$. Hence, by compactness of the latter, there exists a finite open sub-covering, with balls, say, $B_{X_1}, B_{X_2}, \ldots, B_{X_q}$, and corresponding positive values $r_1, r_2, \ldots, r_q$. Let $r_0=\min\left\{r_1, r_2, \ldots, r_q\right\}$, then it follows that for every $Y$ within the $d-$ball centered at $X_0$ and radius $r_0$, i.e., $Y\in B\left(X_0; r_0\right)$, it also holds $f(Y)>f\left(X_0\right)$ and $X_0$ is a point of strict local minimum as asserted. Alternatively, one could also define $r_0$ above as follows: let $r_0(X)=\min\left(1, \max r \text{ such that} f\left(X_0+t\left(X-X_0\right)\right)>f\left(X_0\right), \text{ for } |t|<r\right)$. The latter is a positive, continuous function on (the compact set) $S^{d-1}$. Therefore we have, by way of writing $r_0=:\min\left\{r(X): X\in S^{d-1}\right\}$, that $r_0>0$.

Similarly, in the case labeled as $\mathrm{c}_2)\ a\le b<0$, we have that $f\left(X_0\right)$ is a strict local maximum.

In the next case $\mathrm{c}_3)\ a<0<b$, we choose two points $X_1, X_2\in S^{d-1}$ with $P_k\left(X_1\right)<0$ and $P_k\left(X_2\right)>0$. Then, along the ray emanating from $X_0$ and passing through $X_1$ there exists $r_1>0$ such that $f\left(X_0+t\left(X_1-X_0\right)\right)<f\left(X_0\right)$ for every $t$ with $0<t<r_1$. In a Similar fashion, along the ray emanating from $X_0$ and passing through $X_2$ we have the existence of a positive $r_2>0$ such that $f\left(X_0+t\left(X_2-X_0\right)\right)>f\left(X_0\right)$ for every $t$ with $0<t<r_2$. It necessarily follows that within every neighborhood of $X_0$ there are points $Y, Z$ such that $f(Y)<f\left(X_0\right)<f(Z)$, i.e., $X_0$ is a saddle point as asserted.

We consider next the case $\mathrm{c}_4)\ 0=a<b$ and observe that the even degree polynomial function $P_k$ cannot be of the form $P_k(X)=P_{2p}(X)=\lambda\left(\left(x_1-x_{0_1}\right)^2+\left(x_2-x_{0_2}\right)^2+\cdots+\left(x_d-x_{0_d}\right)^2\right)^p$, by Lemma 6.2. Thus, it is easy to see, first, that the set properly included on the sphere where the $k^{th}$-homogeneous polynomial $P_k$ vanishes cannot be any relative open subset of the sphere either and, second, that such an algebraic set, denoted by $\left(P_k\right)^{-1}(0)\cap S^{d-1}\subset S^{d-1}$, is precisely characterized by the vanishing of the vector function

$$G=\left(G_1, G_2\right): U\to\mathbb{R}^2,$$

where $G_1(X):=\left(x_1-x_{0_1}\right)^2+\left(x_2-x_{0_2}\right)^2+\cdots+\left(x_d-x_{0_d}\right)^2-1$ and $G_2:=P_k$, i.e., by the system

$$\left.\begin{aligned} G_1(X)&:=\left(x_1-x_{0_1}\right)^2+\left(x_2-x_{0_2}\right)^2+\cdots+\left(x_d-x_{0_d}\right)^2-1=0 \\ G_2&:=P_k=P_k(X)=P_k\left(\left(x_1-x_{0_1}\right),\left(x_2-x_{0_2}\right),\ldots,\left(x_d-x_{0_d}\right)\right)=0 \end{aligned}\right\}. \tag{6.4}$$

Then, it is clear that the non-empty solution set of this system can be expressed as a finite union of algebraic submanifolds of the sphere $S^{d-1}$. All of these submanifolds may have diverse dimensions strictly less than $d-1$, the extreme cases being those with dimension $0$, i.e., isolated points, and $d-2$, which occurs if there exist points $X \in G^{-1}(0,0) \subset (P_k)^{-1}(0) \cap S^{d-1}$ where the Jacobian matrix $G'$ attains maximal rank, i.e., $rank(G'(X)) = 2$.

In order to go on with the analysis of this case it is convenient to consider, first, the case where the solution set is constituted only by a finite number of isolated points, say $X_1, X_2, ..., X_q$, and we evaluate at everyone of them the next homogeneous polynomial which appears as non-vanishing in the Taylor´s formula of the function. Let us denote that polynomial by

$$P_l = P_l\left(\left(x_1 - x_{0_1}\right), \left(x_2 - x_{0_2}\right), ..., \left(x_d - x_{0_d}\right)\right) = \sum a_{i_1 i_2 ... i_l}\left(x_{i_1} - x_{0_{i_1}}\right)\left(x_{i_2} - x_{0_{i_2}}\right)...\left(x_{i_l} - x_{0_{il}}\right), \quad (6.5)$$

with $l > k$, $a_{i_1 i_2 ... i_l} = \frac{1}{l!} f_{i_1 i_2 ... i_l}(X_0)$.

Then, it is easy to verify in this particular instance the parts of the statement of the theorem labeled as sub-cases $c_{41})$ through $c_{45})$.

Next, let us consider the case where there exist points $X \in G^{-1}(0,0) = (P_k)^{-1}(0) \cap S^{d-1}$ such that the Jacobian matrix $G'$ attains maximal rank, i.e., $rank(G'(X)) = 2$. As stated previously, we have here the existence of an algebraic submanifold of the solution set of (6.4) with the maximum possible of dimensions, $d-2$. It will be clear, from the exposition that follows, that the other cases of existence of submanifolds with dimensions lesser than that may be treated in a quite similar fashion. Thus, we consider next the

## Second Subsidiary Problem:

Determine the critical points of $P_l = P_l\left(\left(x_1 - x_{0_1}\right), \left(x_2 - x_{0_2}\right), ..., \left(x_d - x_{0_d}\right)\right)$ with the restriction that $X$ belongs to the differentiable manifold $M^{d-2} := \{X : G(X) = 0, rank(G'(X)) = 2\}$, where $G = (G_1, G_2)$, with $G_1(X) := \left(x_1 - x_{0_1}\right)^2 + \left(x_2 - x_{0_2}\right)^2 + \cdots + \left(x_d - x_{0_d}\right)^2 - 1$ and $G_2 := P_k$. It is clear that the above is a $(d-2)-$dimensional algebraic differentiable submanifold contained in the compact, algebraic set $(P_k)^{-1}(0) \cap S^{d-1}$, and that, as in the case of $P_k$, the image of the restriction $P_{l|S^{d-1}}$ is a closed finite interval, according to the statement in Lemma 6.3. Hence, the image of the further restriction $P_{l|M^{d-2}}$ is contained in the latter.

In order to continue the argument in the case $c_4)$ for $P_k$, let us consider, now, the case where there exist points $X \in G^{-1}(0,0) \subset (P_k)^{-1}(0) \cap S^{d-1}$ such that the Jacobian matrix $G'$ does not attain maximal rank, i.e., $rank(G'(X)) < 2$. Since, on the other hand, we are also assuming that the system (6.4) admits a non-empty solution set, it follows that the set satisfying those conditions can be expressed as a finite union of algebraic submanifolds of the sphere $S^{d-1}$ with dimensions strictly less that $d-2$. Therefore, each of those submanifolds may also be represented as the solution set of a system of equations similar to the one described next:

$$\left.\begin{array}{r} G_1(X) := \left(x_1 - x_{0_1}\right)^2 + \left(x_2 - x_{0_2}\right)^2 + \cdots + \left(x_d - x_{0_d}\right)^2 - 1 = 0 \\ G_2 := P_{k_1} = P_{k_1}(X) = P_{k_1}\left(\left(x_1 - x_{0_1}\right), \left(x_2 - x_{0_2}\right), \ldots, \left(x_d - x_{0_d}\right)\right) = 0 \\ \cdots\cdots\cdots\cdots\cdots\cdots\cdots\cdots\cdots\cdots \\ G_{k_r+1} := P_{k_r} = P_{k_r}(X) = P_{k_r}\left(\left(x_1 - x_{0_1}\right), \left(x_2 - x_{0_2}\right), \ldots, \left(x_d - x_{0_d}\right)\right) = 0 \end{array}\right\}. \tag{6.6}$$

Here, all of the polynomial appearing in the left hand side have degrees strictly lower than $k$ and the indicated number of equations, $k_r + 1$, is less than or equal to $d$, $3 \le k_r + 1 \le d$. Besides, the rank of the corresponding Jacobian matrix is maximal, i.e., equal to $k_r + 1$ and, consequently, the previously determined solution set represents an algebraic submanifold of the unit sphere with dimension $d - (k_r + 1)$, $M^{d-(k_r+1)} \subset (P_k)^{-1}(0) \cap S^{d-1} \subset S^{d-1}$.

Moreover we can proceed now, in each of those cases, to consider the corresponding subsidiary problem of determining the critical points of the restricted polynomial $P_{l|_{M^{d-(k_r+1)}}}$, i.e., of the mentioned polynomial $P_l = P_l\left(\left(x_1 - x_{0_1}\right), \left(x_2 - x_{0_2}\right), \ldots, \left(x_d - x_{0_d}\right)\right)$ with the restriction that $X$ belongs to the latter differentiable manifold $M^{d-(k_r+1)}$, following the same procedure as we did before for the case of $M^{d-2}$. It is also clear that the compact, algebraic set $(P_k)^{-1}(0) \cap S^{d-1}$ is equal to the union of those submanifolds together with the previous one $M^{d-2}$.

Thus, with a procedure similar to the previous considerations, regarding the different possibilities for the behavior of $P_k$ on $S^{d-1}$, we may go on to draw conclusions from the behavior of $P_l$ on the compact, algebraic subset $(P_k)^{-1}(0) \cap S^{d-1}$. For example:

1) for $P_l$. If $l$ is odd and there exists $X \in (P_k)^{-1}(0) \cap S^{d-1}$ such that $P_l(X) \neq 0$, then $X_0$ is a saddle point.

2) for $P_l$. If $l$ is even we have again, for the homogeneous polynomial function $P_l$, the five possibilities similar to those described in Lemma 6.2, from which we conclude that:

In case $c_1$) for $P_l$, $f(X_0)$ is a strict local minimum, because we can make the similar kind of reasoning as in the previous case by observing now that, for every $X \in S^{d-1}$, $P_k(X) > 0$ or $P_l(X) > 0$. Hence, by means of a refinement of the argument as in the same case for $P_k(X)$ and by analyzing again the values of the function $f$, one concludes that there exists a positive real number, $r > 0$, such that for every $t$ with $|t| < r$ it holds $f\left(X_0 + t\left(X - X_0\right)\right) > f\left(X_0\right)$. The rest of the argument follows at once.

In case $c_2$) for $P_l$, $f(X_0)$ is of saddle type, since it was previously a candidate for minimum on the basis of the analysis of $P_k$ on $S^d$ but, the further restriction of $P_l$ to $(P_k)^{-1}(0) \cap S^{d-1}$ makes it a prospective in the other way around, i.e., a local maximum.

In case $c_3$) for $P_l$, $X_0$ is obviously a saddle point. And the same can be said of case $c_5$) for $P_l$, because there exists $X \in (P_k)^{-1}(0) \cap S^{d-1}$ such that $P_l(X) \neq 0$.

Now, in case $c_4)$ for $P_l$, $f(X_0)$ is kept as a candidate for furnishing a minimum and, once again, there exists a subset $(P_l)^{-1}(0)\cap(P_k)^{-1}(0)\cap S^{d-1}\subset(P_k)^{-1}(0)\cap S^{d-1}\subset S^{d-1}$, where the $l^{th}$-homogeneous polynomial $P_l$ vanishes, as it was previously the case with $P_k$. In this instance, it is easy to see for the same reasons as before, Lemma 6.2, that the even degree polynomial function $P_l$ neither can be of the form $P_l(X)=P_{2p}(X)=\lambda\left(\left(x_1-x_{0_1}\right)^2+\left(x_2-x_{0_2}\right)^2+\cdots+\left(x_d-x_{0_d}\right)^2\right)^p$ nor be identically vanishing on $(P_k)^{-1}(0)\cap S^{d-1}$. Moreover, $P_l$ cannot vanish on any relative open subset of $(P_k)^{-1}(0)\cap S^{d-1}$.

Observe that, now, the compact, algebraic set $(P_l)^{-1}(0)\cap(P_k)^{-1}(0)\cap S^{d-1}$ is characterized by the vanishing of the vector function

$$G=(G_1,G_2,G_3):U\to\mathbb{R}^3,$$

with $G_1(X):=\left(x_1-x_{0_1}\right)^2+\left(x_2-x_{0_2}\right)^2+\cdots+\left(x_n-x_{0_n}\right)^2-1$, $G_2:=P_k$, $G_3:=P_l$, i.e., by the system

$$\left.\begin{aligned} G_1(X)&:=\left(x_1-x_{0_1}\right)^2+\left(x_2-x_{0_2}\right)^2+\cdots+\left(x_d-x_{0_d}\right)^2-1=0\\ G_2&:=P_k=P_k(X)=P_k\left(\left(x_1-x_{0_1}\right),\left(x_2-x_{0_2}\right),\ldots,\left(x_d-x_{0_d}\right)\right)=0\\ G_3&:=P_l=P_l(X)=P_l\left(\left(x_1-x_{0_1}\right),\left(x_2-x_{0_2}\right),\ldots,\left(x_d-x_{0_d}\right)\right)=0\end{aligned}\right\}. \tag{6.7}$$

Hence, and again as it was previously the case, after equation (6.4), the non-empty solution set of this system can be expressed as a finite union of algebraic submanifolds of $(P_k)^{-1}(0)\cap S^{d-1}\subset S^{d-1}$. All of these submanifolds may have diverse dimensions strictly less than $d-2$, the extreme cases being those with dimension $0$, i.e., isolated points, and $d-3$, which happens if there exist points $X\in G^{-1}(0,0,0)\subset(P_l)^{-1}(0)\cap(P_k)^{-1}(0)\cap S^{d-1}$ where the Jacobian matrix $G'$ attains maximal rank, i.e., for the present case $rank\left(G'(X)\right)=3$.

Finally, it is also easy to repeat the procedure as before, from case 1) for $P_l$ on and, in particular, to verify in this instance the similar parts in the statement of the theorem labeled as $c_{41})$ through $c_{45})$.

It is now quite clear that the above procedure may be repeated again and again, until all possibilities are exhausted after a finite number of steps, since in every one of those further steps the dimensions diminish, thus providing the complete proof for the case labeled at the beginning as $c_4)$.

Analogously, it is also clear that the case $c_5)$ may be treated in a completely similar way to that in the previous case $c_4)$, by simply reversing all of the inequalities; the strong possibility now is that the critical point $X_0$ provides a maximum.

The theorem is proved.

## 7. Final Examples and Exercises.

**Example 7.1.** Let us analyze the following function $f(x,y,z)=5x^2+3y^2+o(2)$, which is already expressed as in Taylor´s formula, where the term $o(2)$ involves, as usual, homogeneous polynomials of degree greater or equal than $3$. Then, the origin with coordinates $(0,0,0)$ is a candidate for extreme and, in the terminology of the previous theorem, $P_2(x,y,z)=5x^2+3y^2$. Thus, we consider first the subsidiary problem of determining the critical point of $P_{2|S^2}:S^2\to\mathbb{R}$, i.e., finding the possible extrema of $P_2$ subject to the restriction $G(x,y,z)=x^2+y^2+z^2-1=0$. Here, the Jacobian matrix may be written $G'(X)=2(x,y,z)$, and we may assume first that $z\neq 0$. Then, in terms of Theorem 3.1 we further obtain successively

$$H(x,y)=(x,y,h(x,y));\ J(x,y)=P_2\circ H(x,y)=5x^2+3y^2$$

$$h_x=-\frac{G_x}{G_z}=-\frac{x}{z};\ h_y=-\frac{G_y}{G_z}=-\frac{y}{z}$$

$$J_x=P_{2,x}+P_{2,z}h_x=10x+0\left(-\frac{x}{z}\right)=10x;\ J_y=P_{2,y}+P_{2,z}h_y=6y+0\left(-\frac{y}{z}\right)=6y$$

Thus, we have to solve the system of polynomial equations

$$\left.\begin{aligned}10x=0\\6y=0\\x^2+y^2+z^2-1=0\end{aligned}\right\}.$$

Either by simple inspection, or by using Swp2.5, we find the solution set of the latter as:

$$\{x=0,y=0,z=1\};\ \{x=0,y=0,z=-1\}.$$

Similarly, for $x\neq 0$, we may write successively

$$H(y,z)=(y,z,h(y,z));\ J(y,z)=P_2\circ H(y,z)=5(h(y,z))^2+3y^2$$

$$h_y=-\frac{G_y}{G_x}=-\frac{y}{x};\ h_z=-\frac{G_z}{G_x}=-\frac{z}{x}$$

$$J_y=P_{2,y}+P_{2,x}h_y=6y+10x\left(-\frac{y}{x}\right)=-4y;\ J_z=P_{2,z}+P_{2,x}h_z=0+10x\left(-\frac{z}{x}\right)=-10z.$$

The system of polynomial equations to solve here is

$$\left.\begin{aligned}-4y=0\\-10z=0\\x^2+y^2+z^2-1=0\end{aligned}\right\},$$

with solution set $\{x=1,y=0,z=0\}$; $\{x=-1,y=0,z=0\}$.

Finally, for $y\neq 0$, one obtains

$$H(x,z)=(x,z,h(x,z));\ J(x,z)=P_2\circ H(x,z)=5x^2+3(h(x,z))^2$$

$$h_x = -\frac{G_x}{G_y} = -\frac{x}{y}\,;\; h_z = -\frac{G_z}{G_y} = -\frac{z}{y}$$

$$J_x = P_{2,x} + P_{2,y}h_x = 10x + 6y\left(-\frac{x}{y}\right) = 4x\,;\; J_z = P_{2,z} + P_{2,y}h_z = 0 + 6y\left(-\frac{z}{y}\right) = -6z$$

The system of polynomial equations to solve now is

$$\left.\begin{array}{c} 4x = 0 \\ -6z = 0 \\ x^2 + y^2 + z^2 - 1 = 0 \end{array}\right\},$$

with solution set $\{x = 0, y = 1, z = 0\}$; $\{x = 0, y = -1, z = 0\}$.

Then, evaluating the polynomial $P_2(x, y, z) = 5x^2 + 3y^2$ at the points obtained gives:

$$P_2(0,0,1) = 0\,, P_2(0,0,-1) = 0\,, P_2(1,0,0) = 5,$$
$$P_2(-1,0,0) = 5,\; P_2(0,1,0) = 3\; P_2(0,-1,0) = 3,$$

and it follows that $P_2(S^2) = [0,5]$.

Hence, we are in the situation described in Theorem 6.5 as case $c_4)$, the origin is a candidate for a minimum but we need to analyze higher degree terms at the points on the 2-sphere where $P_2$ vanishes. In the present case there are only two points as described above $(0,0,1)$ and $(0,0,-1)$. Thus, we consider next

**Example 7.1.1.** Let us modify the previous example by assuming now that the function is expressed by $f(x, y, z) = 5x^2 + 3y^2 + x^2y + 2xz^2 + z^3 + o(3)$. Then it is easy to see that the homogeneous third degree polynomial $P_3(x, y, z) = x^2y + 2xz^2 + z^3$, attains the values $1$ and $-1$ respectively at the points where the previous one $P_2$ was vanishing. Hence, the function in this particular example has a saddle point at the origin, since we are in the situation described in Theorem 6.5, case $c_{42})$.

**Example 7.1.2.** Let us next modify the previous example, by also adding higher order terms, as follows $f(x, y, z) = 5x^2 + 3y^2 + x^2y + 2xz^2 + 5x^3y + 7xz^3 + 3z^4 + o(4)$. Then, we see that the third degree homogeneous polynomial $P_3(x, y, z) = x^2y + 2xz^2$ is vanishing at the points $(0,0,1)$ and $(0,0,-1)$, which allows to keep the expectation that the origin may furnish a minimum for the function and, in fact, by analyzing the values of the next, fourth order homogeneous polynomial, i.e., $P_4(x, y, z) = 5x^3y + 7xz^3 + 3z^4$ we find that $P_4(0,0,1) = P_4(0,0,-1) = 3$. Thus, this belongs to the situation described in Theorem 6.5, cases $c_{41})$, second part, and $c_{43})$. The function has a strict minimum at the origin.

**Example 7.1.3.** By the same token if we have, for the given function, situations like

**7.1.3.1.** $f(x,y,z)=5x^2+3y^2+g(x,y,z)$, where $g(x,y,z)$ is a (not necessarily finite) sum of terms containing, everyone of them, at least one of the variables $x$ or $y$ as a factor, then, since $g(0,0,1)=g(0,0,-1)=0$, we are in the situation described in Theorem 6.5, case $\mathrm{c}_{41})$, and the origin is a point of non-strict minimum for the function. In particular, $f(0,0,z)=0$ for every $z\in\mathbb{R}$.

**7.1.3.2.** $f(x,y,z)=5x^2+3y^2+g(x,y,z)+z^{2p}+\text{higher degree terms}$, where $g(x,y,z)$ is a sum of terms containing, everyone of them, at least one of the variables $x$ or $y$, with degree strictly greater than $2$ and less than or equal to $2p$, $p$ a natural number, $p\in\mathbb{N}, p>1$, then the origin is a point of strict minimum for the function. Theorem 6.5, cases $\mathrm{c}_{41})$ and $\mathrm{c}_{43})$.

**7.1.3.3.** $f(x,y,z)=5x^2+3y^2+g(x,y,z)+z^{2p+1}+\text{higher degree terms}$, where $g(x,y,z)$ is a sum of terms containing, everyone of them, at least one of the variables $x$ or $y$, here with degree strictly greater than $2$ and less than or equal to $2p+1$, $p$ a natural number, $p\in\mathbb{N}, p>1$, then the origin is a saddle point for the function. Theorem 6.5, cases $\mathrm{c}_{41})$ and $\mathrm{c}_{42})$.

**Exercise 7.2.** Prove that every function of the form $f(x,y,z)=5x^2+3y^2+7z^2+o(2)$ has a strict minimum at the origin of coordinates, $(0,0,0)$, by showing that in these cases the image under the homogeneous polynomial $P_2(x,y,z)=5x^2+3y^2+7z^2$ of the 2-sphere, $S^2$, is the closed interval $[3,7]$. Theorem 6.5, case $\mathrm{c}_1)$.

**Example 7.3.** Let us consider now the function $f(x,y,z)=x^2+2xy+y^2+o(2)$, where the first homogeneous polynomial is $P_2(x,y,z)=x^2+2xy+y^2$. Again, the subsidiary problem in this case consists in determining the critical point of $P_{2|S^2}:S^2\to\mathbb{R}$, i.e., finding the possible extrema of $P_2$ subject to the restriction $G(x,y,z)=x^2+y^2+z^2-1=0$, whose Jacobian matrix is given by $G'(X)=2(x,y,z)$, and we assume first that $z\neq 0$. Then, in terms of Theorem 3.1 we further obtain

$$H(x,y)=(x,y,h(x,y));\; J(x,y)=P_2\circ H(x,y)=x^2+2xy+y^2$$

$$h_x=-\frac{G_x}{G_z}=-\frac{x}{z};\; h_y=-\frac{G_y}{G_z}=-\frac{y}{z}$$

$$J_x=P_{2,x}+P_{2,z}h_x=2x+2y+0\left(-\frac{x}{z}\right)=2x+2y;$$

$$J_y=P_{2,y}+P_{2,z}h_y=2x+2y+0\left(-\frac{y}{z}\right)=2x+2y$$

Thus, we have to solve the system of polynomial equations

$$\left.\begin{array}{r}2x+2y=0\\2x+2y=0\\x^2+y^2+z^2-1=0\end{array}\right\},$$

whose solution set we find either by simple inspection, or by using Swp2.5, to be :

$$\left\{x=\frac{1}{2}\sqrt{-2z^2+2}, y=-\frac{1}{2}\sqrt{-2z^2+2}, z=z\right\};$$
$$\left\{x=-\frac{1}{2}\sqrt{-2z^2+2}, y=\frac{1}{2}\sqrt{-2z^2+2}, z=z\right\}.$$

Obviously, this representation of the solution, which is even valid or extendable by continuity to $z=0$, and/or also included in the next case considered below ( $x\neq 0$ ), provides us with, and may be interpreted as, a parametrization of the 1-dimensional manifold, $M^1 := S^2 \cap \{(x,y,z): x+y=0\}$, obtained by intersecting the 2-sphere $S^2$ with the plane determined by equation $x+y=0$. The parameter $z$ taking values in the interval $[-1,1]$.

Similarly, for $x\neq 0$, we may write

$$H(y,z)=(y,z,h(y,z));\ J(y,z)=P_2\circ H(y,z)=(h(y,z))^2+2(h(y,z))y+y^2$$
$$h_y=-\frac{G_y}{G_x}=-\frac{y}{x};\ h_z=-\frac{G_z}{G_x}=-\frac{z}{x}$$
$$J_y=P_{2,y}+P_{2,x}h_y=2x+2y+(2x+2y)\left(-\frac{y}{x}\right)=2x-2\frac{y^2}{x};$$
$$J_z=P_{2,z}+P_{2,x}h_z=0+(2x+2y)\left(-\frac{z}{x}\right)=-2z-2\frac{yz}{x}.$$

Now, the system of equations to solve is

$$\left.\begin{array}{r} 2x-2\dfrac{y^2}{x}=0 \\ -2z-2\dfrac{yz}{x}=0 \\ x^2+y^2+z^2-1=0 \end{array}\right\},$$

with solution set:

$$\left\{x=\frac{1}{2}\sqrt{-2z^2+2}, y=-\frac{1}{2}\sqrt{-2z^2+2}, z=z\right\};$$
$$\left\{x=-\frac{1}{2}\sqrt{-2z^2+2}, y=\frac{1}{2}\sqrt{-2z^2+2}, z=z\right\};$$
$$\left\{x=\frac{1}{2}\sqrt{2}, y=\frac{1}{2}\sqrt{2}, z=0\right\};\ \left\{x=-\frac{1}{2}\sqrt{2}, y=-\frac{1}{2}\sqrt{2}, z=0\right\}.$$

Finally, for $y\neq 0$, we may proceed to perform again the computations as above or, otherwise, observe that this problem is symmetric in the variables $x, y$ so that the solution set, obtained after repeating the procedure, is the same as the latter one.

Next, we evaluate the polynomial $P_2(x,y,z) = x^2 + 2y + y^2$ at the solution set, obtaining: for points on $M^1$, $P_2(x,y,z) = x^2 + 2y + y^2 = (x+y)^2 = 0$, while for both of the remaining two points we get $P_2(x,y,z) = x^2 + 2y + y^2 = \left(\sqrt{2}\right)^2 = 2$.

It follows that the image of the 2-sphere $S^2$ under $P_2$ is $P_2(S^2) = [0,2]$, i.e., the close interval $[0,2]$. Hence, we are in the situation described in Theorem 6.5 as case $c_4$), the origin is a candidate for a minimum but we need to analyze higher degree terms at the points on the 2-sphere where $P_2$ vanishes.

**Example 7.3.1.** Let us modify the latter example by assuming now that the function to be considered is $f(x,y,z) = x^2 + 2xy + y^2 + x^4 + y^4 + o(4)$. Thus, we can use the information obtained previously and proceed to study the behavior of the homogeneous polynomial $P_4(x,y,z) = x^4 + y^4$ on the subset of the 2-sphere $S^2$ where $P_2$ vanishes, i.e., in the previous notation the compact, algebraic 1-dimensional submanifold of the 2-sphere, $P_2^{-1}(0) \cap S^2 = M^1$. We may do so by considering the subsidiary problem of finding the critical points of $P_4(x,y,z) = x^4 + y^4$ subject to the restrictions determined by the simultaneous equations: $G_1(x,y,z) = G(x,y,z) = x^2 + y^2 + z^2 - 1 = 0$, $G_2(x,y,z) = x + y = 0$. We leave it as an exercise to show that the solution set of this problem is constituted by the four points:

$$\left(\frac{1}{2}\sqrt{2}, -\frac{1}{2}\sqrt{2}, 0\right);\ \left(-\frac{1}{2}\sqrt{2}, \frac{1}{2}\sqrt{2}, 0\right);\ (0,0,1);\ (0,0,-1);$$

the values of $P_4$ at those points being:

$$P_4\left(\frac{1}{2}\sqrt{2}, -\frac{1}{2}\sqrt{2}, 0\right) = P_4\left(-\frac{1}{2}\sqrt{2}, \frac{1}{2}\sqrt{2}, 0\right) = \frac{1}{2};\ P_4(0,0,1) = P_4(0,0,-1) = 0.$$

Then it follows that the image of the compact, algebraic submanifold $M^1$ under $P_4$ is: $P_4(M^1) = [0, 1/2]$. Thus, in order to conclude whether the origin is a minimum, or not, we would have to analyze higher degree terms of the function, but now only at the two points where $P_4$ is vanishing. So that this situation is quite similar to the previous one in exercises 7.1 and followings, and it is left as an exercise for the interested reader the construction of the further, various possibilities.

**Example 7.4.** We increase now the ambient space dimension and consider the function initially given by $f(w,x,y,z) = x^2 + 2xy + y^2 + o(2)$. Homogeneous polynomial $P_2(w,x,y,z) = x^2 + 2xy + y^2$. Subsidiary problem: determine the critical point of $P_{2|S^3} : S^3 \to \mathbb{R}$, i.e., find the critical points of $P_2$ subject to the restriction $G(w,x,y,z) = w^2 + x^2 + y^2 + z^2 - 1 = 0$, whose Jacobian matrix is given by $G'(X) = 2(w,x,y,z)$, and we analyze first the region where $z \neq 0$. By Theorem 3.1 we write

$$H(w,x,y) = (w,x,y,h(w,x,y));\ J(w,x,y) = P_2 \circ H(w,x,y) = x^2 + 2xy + y^2$$

$$h_w = -\frac{G_w}{G_z} = -\frac{w}{z};\ h_x = -\frac{G_x}{G_z} = -\frac{x}{z};\ h_y = -\frac{G_y}{G_z} = -\frac{y}{z}$$

$$J_w = P_{2,w} + P_{2,w}h_w = 0 + 0\left(-\frac{w}{z}\right) = 0;$$

$$J_x = P_{2,x} + P_{2,z}h_x = 2x + 2y + 0\left(-\frac{x}{z}\right) = 2x + 2y;$$

$$J_y = P_{2,y} + P_{2,z}h_y = 2x + 2y + 0\left(-\frac{y}{z}\right) = 2x + 2y.$$

Thus, we have to solve the system of polynomial equations

$$\left.\begin{array}{r} 0 = 0 \\ 2x + 2y = 0 \\ 2x + 2y = 0 \\ w^2 + x^2 + y^2 + z^2 - 1 = 0 \end{array}\right\}.$$

The solution set here is:

$$\left.\begin{array}{l} \left\{ w = w, x = \frac{1}{2}\sqrt{-2z^2 - 2w^2 + 2}, y = -\frac{1}{2}\sqrt{-2z^2 - 2w^2 + 2}, z = z \right\} \\ \left\{ w = w, x = -\frac{1}{2}\sqrt{-2z^2 - 2w^2 + 2}, y = \frac{1}{2}\sqrt{-2z^2 - 2w^2 + 2}, z = z \right\} \end{array}\right\} \quad (7.1)$$

It is rather obvious, for reasons of symmetry, that we get the same result at the region where $w \neq 0$. On the other hand, when applied $P_2(w,x,y,z) = x^2 + 2xy + y^2$ to the elements of this solution one obtains, obviously, $P_2(w,x,y,z) = (x+y)^2 = 0$.

Let us analyze next the regions where $x \neq 0$:

$$H(w,y,z) = (w, h(y,z), y, z);$$

$$J(w,y,z) = P_2 \circ H(w,y,z) = (h(w,y,z))^2 + 2(h(w,y,z))y + y^2$$

$$h_w = -\frac{G_w}{G_x} = -\frac{w}{x};\ h_y = -\frac{G_y}{G_x} = -\frac{y}{x};\ h_z = -\frac{G_z}{G_x} = -\frac{z}{x}$$

$$J_w = P_{2,w} + P_{2,x}h_w = 0 + (2x + 2y)\left(-\frac{w}{x}\right) = -2w - 2\frac{yw}{x}$$

$$J_y = P_{2,y} + P_{2,x}h_y = 2x + 2y + (2x + 2y)\left(-\frac{y}{x}\right) = 2x - 2\frac{y^2}{x};$$

$$J_z = P_{2,z} + P_{2,x}h_z = 0 + (2x + 2y)\left(-\frac{z}{x}\right) = -2z - 2\frac{yz}{x}$$

Therefore, the system of equations to solve is

$$\left.\begin{array}{c} -2w-2\dfrac{yw}{x}=0 \\ 2x-2\dfrac{y^2}{x}=0 \\ -2z-2\dfrac{yz}{x}=0 \\ w^2+x^2+y^2+z^2-1=0 \end{array}\right\},$$

whose solution set is the union of the ones described in the following equations:

$$\left.\begin{array}{c} \left\{w=\sqrt{-2y^2-z^2+1}, x=-y, y=y, z=z\right\} \\ \left\{w=-\sqrt{-2y^2-z^2+1}, x=-y, y=y, z=z\right\} \end{array}\right\} \tag{7.2}$$

$$\left.\begin{array}{c} \left\{w=0, x=\frac{1}{2}\sqrt{-2z^2+2}, y=-\frac{1}{2}\sqrt{-2z^2+2}, z=z\right\} \\ \left\{w=0, x=-\frac{1}{2}\sqrt{-2z^2+2}, y=\frac{1}{2}\sqrt{-2z^2+2}, z=z\right\} \end{array}\right\} \tag{7.3}$$

$$\left.\begin{array}{c} \left\{w=0, x=\frac{1}{2}\sqrt{2}, y=\frac{1}{2}\sqrt{2}, z=0\right\} \\ \left\{w=0, x=-\frac{1}{2}\sqrt{2}, y=-\frac{1}{2}\sqrt{2}, z=0\right\} \end{array}\right\} \tag{7.4}$$

$$\left.\begin{array}{c} \left\{w=0, x=\frac{1}{2}\sqrt{2}, y=-\frac{1}{2}\sqrt{2}, z=0\right\} \\ \left\{w=0, x=-\frac{1}{2}\sqrt{2}, y=\frac{1}{2}\sqrt{2}, z=0\right\} \end{array}\right\} \tag{7.5}$$

If we analyze now the regions where $y \neq 0$ we obtain, for reasons of symmetry of the four variables involved, by interchanging $x \leftrightarrow y$ and $w \leftrightarrow z$ the following additional solution sets:

$$\left.\begin{array}{c} \left\{w=w, x=x, y=-x, z=\sqrt{-2x^2-w^2+1}\right\} \\ \left\{w=w, x=x, y=-x, z=-\sqrt{-2x^2-w^2+1}\right\} \end{array}\right\} \tag{7.6}$$

$$\left.\begin{array}{c} \left\{w=w, x=\frac{1}{2}\sqrt{-2w^2+2}, y=-\frac{1}{2}\sqrt{-2w^2+2}, z=0\right\} \\ \left\{w=w, x=-\frac{1}{2}\sqrt{-2w^2+2}, y=\frac{1}{2}\sqrt{-2w^2+2}, z=0\right\} \end{array}\right\} \tag{7.7}$$

It is easy to see that the polynomial $P_2(w,x,y,z) = x^2 + 2xy + y^2 = (x+y)^2$ vanishes at almost all subsets above, except at the points $\left(0,-\sqrt{2}/2,-\sqrt{2}/2,0\right)$ and $\left(0,\sqrt{2}/2,\sqrt{2}/2,0\right)$, described in (7.4), where it takes the value $P_2\left(0,-\sqrt{2}/2,-\sqrt{2}/2,0\right) = P_2\left(0,\sqrt{2}/2,\sqrt{2}/2,0\right) = 2$. Therefore, the image of $S^3$ under $P_2$ is the close interval $P_2\left(S^3\right) = [0,2]$ and, consequently, the origin $(0,0,0,0)$ is a candidate for a minimum, but one needs to analyze higher order homogeneous polynomials at the points of the sphere $S^3$ where $P_2$ vanishes, i.e., at the set $S^3 \cap P_2^{-1}(0)$, which is easily seen to be the union of a compact, 2-dimensional algebraic submanifold of $S^3$, say $M^2$ described in (7.1), (7.2) and (7.6) above, the union of two 1-dimensional ones, $M^1$, described in (7.3) and (7.7), and two isolated points, $\left(0,\sqrt{2}/2,-\sqrt{2}/2,0\right)$ and $\left(0,-\sqrt{2}/2,\sqrt{2}/2,0\right)$ given in (7.5). However observe, too, that those two points and $M^1$ are subsets of $M^2$, so that for analyzing further higher degree homogeneous polynomials it is enough to consider only the latter set.

**Exercise 7.4.1.1.** Prove that the functions of the general form

$$f(w,x,y,z) = x^2 + 2xy + y^2 + z^{2p+1} + o(2p+1),$$

$p \in \mathbb{N}, p \geq 1$, have a saddle point at the origin of coordinates $(0,0,0,0)$. Theorem 6.5, case $c_{42})$.

**Exercise 7.4.1.2.** Prove that the functions of the general form

$$f(w,x,y,z) = x^2 + 2xy + y^2 + z^{2p+2} + o(2p+2),$$

with $p \in \mathbb{N}, p \geq 1$, do not necessarily have a point or strict minimum at the origin of coordinates $(0,0,0,0)$. Hint: see carefully the situation described in Theorem 6.5, case $c_{43})$.

**Exercise 7.4.1.3.** Prove that the functions of the general form

$$f(w,x,y,z) = x^2 + 2xy + y^2 + g(w,x,y,z),$$

where $g(w,x,y,z)$ is a sum of homogeneous polynomials of degrees greater than two, containing everyone of them the monomial $(x+y)$ as a multiplicative factor, fulfills the situation described in Theorem 6.5, case $c_{41})$, so that the origin is a point of non-strict minimum for those functions and, in particular, $f(w,x,y,z) = 0$ for every $(w,x,y,z)$ such that $x+y=0$. Observation: the sum indicated for $g(w,x,y,z)$ could contain only a finite number of terms, i.e., a polynomial function, or else constituted in part by a convergent power series, for example $g(w,x,y,z) = wz\sum_{n=1}^{\infty}\frac{1}{n!}(x+y)^n$.

**Exercise 7.4.1.4.** Prove that the function $f(w,x,y,z) = x^2 + 2xy + y^2 + xz^3 + z^4 + o(4)$ exhibits a saddle point at the origin. Hint: Theorem 6.5, case $c_{44})$,

**Example 7.4.2.** Let us assume next that $f(w,x,y,z) = x^2 + 2xy + y^2 + z^4\left(x^2+y^2\right) + o(6)$.

Here the homogeneous polynomial to analyze is $P_6(w,x,y,z)=z^4(x^2+y^2)$. In this case, in view of the result in Example 7.4, mainly the observation at the end, we have to consider the following subsidiary problem: to determine the critical point of $P_{6|M^2}:M^2\to\mathbb{R}$.

It is easy to see, by equations (7.1), (7.2) and (7.3) above that $M^2$ is given by $M^2=\{(w,x,y,z):w^2+x^2+y^2+z^2-1=0,x+y=0\}$, so that we consider the corresponding problem of determining the critical points of the homogeneous polynomial $P_6$ when restricted by the equations $G_1(x,y,z)=x^2+y^2+z^2-1=0$ and $G_2(x,y,z)=x+y=0$.

We leave it as an exercise for the reader the calculation of the above subsidiary problem, by using Theorem 3.1, and observe that there exist points on $M^2$ where $P_6$ is also vanishing, so that in order to follow up with the analysis of the case one would need to consider even higher order homogeneous polynomials, provided these are known.

**Exercise 7.4.3.** Let $f(w,x,y,z)=x^2+2xy+y^2+z^4(x^2+y^2)+z^4w^4+w^8+o(8)$. Follow up the analysis of the case, by considering now the homogeneous polynomial $P_8(w,x,y,z)=z^4w^4+w^8$, restricted to the set $P_6^{-1}(0)\cap M^2$, by means of the corresponding subsidiary problems as needed.

**Exercise 7.5.** Consider the two similar problems $f(w,x,y,z)=x^4+y^4+z^6+x^2z^4+w^8+o(8)$ and $f(w,x,y,z)=x^4+y^4+z^6+x^2z^4+w^9+w^3x^3z^3+o(9)$. Determine, in each case, if the origin is a maximum, minimum or saddle point.

Hint: in the second step of both cases you have to determine the points where the polynomial functions $G_1(w,x,y,z)=w^2+x^2+y^2+z^2-1$ and $G_2(w,x,y,z)=x^4+y^4$ are simultaneously vanishing. Let $G=(G_1,G_2)$ and observe that there exist no points $X:=(w,x,y,z)$ that fulfill both conditions $G(X)=0$ and $rank(G'(X))=2$. In fact observe too that, when restricted to the real field, $x^4+y^4=0$ if, and only if, $x=0$ and $y=0$, so that the solution set of $G_1=0,G_2=0$ is the same as that of $G_1(w,x,y,z)=w^2+x^2+y^2+z^2-1=0,x=0,y=0$. Verify that the rank of the Jacobian associated to the latter system of equations is maximal, equal to $3$, so that it defines a one-dimensional compact, algebraic submanifold of the unit sphere. Therefore, in the following step you have to study the behavior of the next nonvanishing homogeneous polynomial restricted to that submanifold.

**Exercise 7.6.** Consider the function $F(x_1,x_2)=x_2^4-x_1^4-x_2^8+x_1^{10}$ exhibited as Example 1 in the article by E. Constantin [4]. By observing that, in our terminology, $P_4(x_1,x_2)=x_2^4-x_1^4$, prove that the origin $(0,0)$ is a saddle point, since $P_4(1,0)=-1$ and $P_4(0,1)=1$, case $c_3$) in Theorem 6.5, without regards to what the terms with degrees higher than $4$ may be, (*cf.* [4], pp. 45-46).

**Exercise 7.7.** Consider the following related problems

7.7.1.
$$\begin{aligned}f(t,u,v,w,x,y,z)&=t^2v^2+t^2w^2+t^2x^2+t^2y^2+t^2z^2+u^2v^2+u^2w^2+u^2x^2+u^2y^2+\\&+u^2z^2+v^4+2v^2w^2+v^2x^2+v^2y^2+v^2z^2+w^4+w^2x^2+w^2y^2+\\&+w^2z^2+7x^6+4y^6+2t^2u^4+x^6y+8xz^6-6z^7+o(7)\end{aligned}$$

7.7.2.
$$f(t,u,v,w,x,y,z) = t^2v^2 + t^2w^2 + t^2x^2 + t^2y^2 + t^2z^2 + u^2v^2 + u^2w^2 + u^2x^2 + u^2y^2 + \\ + u^2z^2 + v^4 + 2v^2w^2 + v^2x^2 + v^2y^2 + v^2z^2 + w^4 + w^2x^2 + w^2y^2 + \\ + w^2z^2 + 7x^6 + 4y^6 + 2t^2u^4 + x^6y + 8xz^6 + 6z^{10} + o(10)$$

7.7.3.
$$f(t,u,v,w,x,y,z) = t^2v^2 + t^2w^2 + t^2x^2 + t^2y^2 + t^2z^2 + u^2v^2 + u^2w^2 + u^2x^2 + u^2y^2 + \\ + u^2z^2 + v^4 + 2v^2w^2 + v^2x^2 + v^2y^2 + v^2z^2 + w^4 + w^2x^2 + w^2y^2 + \\ + w^2z^2 - 2t^2u^2v^2 + vz^6 + 3wx^3y^3 + x^6y - 8xz^6 + 5u^{10} + o(10)$$

Determine, in each case, if the origin is a maximum, minimum or saddle point.

Departamento de Matemáticas, Facultad. de Ciencias Exactas, Físicas y Naturales – Universidad Nacional de Córdoba, Ciudad Universitaria, Córdoba, ARGENTINA – e-mail: sgigena@efn.uncor.edu

Departamento de Matemáticas, Facultad de Ciencias Exactas, Ingeniería y Agrimensura – Universidad Nacional de Rosario, Rosario, ARGENTINA – e-mail: sgigena@fceia.unr.edu.ar